\documentclass[a4paper,10pt]{amsart}
\usepackage{amssymb}
\usepackage{booktabs}
\usepackage{hyperref}
\newtheorem{theorem}{Theorem}[section]
\newtheorem{proposition}[theorem]{Proposition}
\newtheorem{lemma}[theorem]{Lemma}
\newtheorem{corollary}[theorem]{Corollary}
\theoremstyle{definition}
\newtheorem{definition}[theorem]{Definition}
\newtheorem*{assumption}{Standing assumption}
\theoremstyle{remark}

\newtheorem{example}[theorem]{Example}
\allowdisplaybreaks[4]

\newcommand{\sgn}{\operatorname{sgn}}

\title{Fourier Positivity and Sharp Poincar\'e Inequalities on Graphs}
\author{Bang-Xian Han}
\address[Bang-Xian Han]{hanbx@sdu.edu.cn,  School of Mathematics, Shandong University, China}
\author{Yuan-He Wang}
\address[Yuan-He Wang]{yuanhewang@mail.sdu.edu.cn,  School of Mathematics, Shandong University, China}
\hypersetup{pdfauthor={Bang-Xian Han and Yuan-He Wang}}
\date{\today}

\begin{document}

\begin{abstract}
We prove sharp dimension-free $L^1$-Poincar\'e inequalities on Hamming cubes, a sharp quadratic $L^1$-Poincar\'e inequality on product Poisson spaces, and a sharp pseudo-Poincar\'e inequality on finite abelian Cayley graphs.
These inequalities use a generalization of the Beckmann boundary introduced by Ivanisvili, Xie, and Zhang~\cite{IvanisviliXieZhang2026}.
In particular, we prove that the optimal dimension-free constant in the classical $L^1$-Poincar\'e inequality on the Hamming cube is $\sqrt{\pi/2}$, answering a question of Ivanisvili et al.~\cite{IvanisviliLiVanHandelVolberg2018}.

Semigroup duality and explicit Malliavin weights reduce these estimates to integrals of first absolute moments of centered Bernoulli and Poisson sums.
The key ingredient is a Fourier positivity principle, which shows that replacing the coefficients by their absolute values can only increase these integrals.
For the underlying Bernoulli functional at subquadratic exponents, we also characterize the maximizers and determine their exact phase transitions.
\end{abstract}

\maketitle

\setcounter{tocdepth}{1}
\tableofcontents

\section{Introduction}

\subsection{Background and main results}
We study $L^1$-type Beckmann Poincaré and pseudo-Poincaré
inequalities on weighted graphs, with particular emphasis on sharp
constants that hold uniformly over several families of graphs.

Let \((V,b,\mu)\) represent a connected, locally finite weighted graph.
Write \(\Delta\) for the Neumann realization of its graph Laplacian
and \(P_t=e^{t\Delta}\) for the associated heat semigroup. For \(0<T\le\infty\), we consider $L^1$-type inequalities of the form
\[
\|g-P_Tg\|_{L^1(\mu)}
\le C\,\mathsf B_q(g),
\]
where \(\mathsf B_q(g)\) is the $q$-Beckmann boundary of \(g\). The case \(T=\infty\) gives a Beckmann Poincaré inequality, while \(0<T<\infty\) corresponds to a Beckmann pseudo-Poincaré inequality.

The Beckmann boundary is a variational relaxation of the usual
gradient boundary. Roughly speaking, it minimizes the \(L^1\)-cost
of \(|v|_q\) over edge vector fields \(v\) satisfying
\(\operatorname{div}v=\Delta g\).
The underlying variational principle goes back to Beckmann's
formulation of optimal transport~\cite{Beckmann1952}.
In the continuous setting, prescribed-divergence flux minimization is related to optimal transport with distance cost; see~\cite[Section~4.2]{Santambrogio2015} for its vector-measure formulation and transport duality.
Motivated by this formulation, Ivanisvili, Xie, and Zhang
\cite{IvanisviliXieZhang2026} introduced such boundary with prescribed divergence $\Delta g$
on the Hamming cube. We extend their definition to weighted graphs.

The Beckmann boundary is useful for two reasons. First, the gradient
field is admissible, so every Beckmann inequality implies the
corresponding classical \(L^1\)-Poincaré or pseudo-Poincaré
inequality with the same constant. Second, our duality theorem
identifies the optimal constant in a Beckmann inequality with that
in a uniform \(L^\infty\)-gradient estimate for the integrated heat
semigroup. The boundary is therefore the smallest term obtained by
optimizing over vector fields in this semigroup-duality argument;
see~\cite{Ledoux2003,Wu2009,DjelloutWu2011,Liu2024} for related
semigroup and Poisson-equation methods.

\medskip
\noindent\textbf{Hamming cubes.}
The Hamming cube with uniform measure is our principal motivating example. Ben Efraim and Lust-Piquard~\cite{BenEfraimLustPiquard2008} proved the dimension-free inequality
\[
\|g-\mu(g)\|_{L^1(\mu)}
\le
\frac{\pi}{2}\,\mathbb E|\nabla_2g|
\]
for all real-valued functions \(g\). Ivanisvili et al.~\cite{IvanisviliLiVanHandelVolberg2018}
improved the upper bound to a value strictly below \(\pi/2\) and
asked whether the sharp Gaussian constant \(\sqrt{\pi/2}\) holds
uniformly over the dimension. Ivanisvili and Stone~
\cite{IvanisviliStone2024} later obtained a further quantitative
improvement.

For Boolean-valued functions, Durcik et al.~
\cite{DurcikIvanisviliRoosXie2026} recently proved that the optimal
constant in the classical inequality is \(1\).
In the Beckmann setting, Ivanisvili, Xie, and Zhang~
\cite[Theorem~3.2]{IvanisviliXieZhang2026} obtained a finite-time
quadratic estimate using the heat semigroup. Taking \(T\to\infty\)
and translating their result into our normalization gives
\[
\|g-\mu(g)\|_{L^1(\mu)}
\le
\frac{\pi}{2}\,\mathsf B_2(g)
\]
for all Boolean-valued functions \(g\).
We establish the following sharp Beckmann Poincaré inequalities
for arbitrary real-valued functions.

\begin{theorem}[Hamming cubes, Corollary~\ref{cor:hamming-cubes}]
Let \(\mathsf B_q\) be the \(q\)-Beckmann boundary from Definition~\ref{def:beckmann-boundary} specialized to the Hamming cube $\{-1,1\}^N$ with the uniform probability measure, and let \(C_q\) denote the optimal constant determined in Theorems~\ref{thm:bernoulli-l2} and~\ref{thm:exact-phase-transitions}.
For \(1\le q\le2\), the Beckmann Poincaré inequality
\[
\|g-\mu(g)\|_{L^1(\mu)}
\le
2^{2/q-1}C_q\,\mathsf B_q(g),\qquad g:\{-1,1\}^N\to \mathbb R
\]
holds on all such spaces.
Moreover, the constant \(2^{2/q-1}C_q\) is optimal uniformly over \(N\).
\end{theorem}

More explicitly, the constants are given by
\[
C_q=\sup_{N\ge1}\frac{S_N}{N^{1/q}},\qquad
S_N=\sum_{k=0}^{N-1}\frac{(N-1)!}{(N-1-k)!N^k},\qquad
C_2=\sqrt{\frac\pi2}.
\]
Since $S_1=1$ and $S_2=3/2$, the first transition occurs at
$q_1=1/\log_2(3/2)$, and $C_q=1$ for $1\le q\le q_1$.

Choosing the admissible vector fields $v=\nabla f$ and $v=(2\nabla f)^+$ in the
definition of $\mathsf B_q(f)$ gives, respectively,
\begin{gather*}
\mathbb E|f-\mu(f)|
\le \frac{C_q}{2}
\mathbb E
\left(\sum_{i=1}^{N}|f(x^{\oplus i})-f(x)|^q\right)^{1/q},\\
\mathbb E|f-\mu(f)|
\le C_q
\mathbb E 
\left(\sum_{i=1}^{N}|(f(x^{\oplus i})-f(x))^+|^q\right)^{1/q}.
\end{gather*}
For $q=2$, the first inequality resolves the question posed in~\cite{IvanisviliLiVanHandelVolberg2018}.
Both displayed constants are sharp for $1\le q\le q_1$, since an indicator depending on one coordinate gives equality.
At $q=q_1$, taking $f=-\mathbf1_A$ in the second inequality gives
\[
\mathbb E h_A^{\log_2(3/2)}\ge2\mu(A)(1-\mu(A)),\qquad
h_A(x)=\mathbf1_A(x)\#\{y\sim x:y\notin A\}.
\]
This recovers the isoperimetric inequality of Kahn and Park~\cite[Theorem~1.1]{KahnPark2020}, while the preceding gradient inequality provides its real-valued functional extension.
Its critical exponent $1/q_1=\log_2(3/2)$ thus corresponds to the first transition between extremizers with one and two nonzero coefficients.

\medskip
\noindent\textbf{Product Poisson spaces.}
Bobkov and Götze~\cite{BobkovGotze1999} studied discrete
isoperimetric and Poincaré-type inequalities for product measures.
Their Proposition~3.6 gives a Gaussian-profile inequality for the forward boundary.
For \(\mu=\operatorname{Pois}(1)^{\otimes N}\), it implies, in our normalization,
\[
\|g-\mu(g)\|_{L^1(\mu)}
\le
\sqrt{2\pi}\mathbb E|\nabla_2g|
\]
for all Boolean-valued functions \(g\).

\begin{theorem}[Product Poisson spaces, Corollary~\ref{cor:product-poisson-spaces}]
Let \(\mathsf B_q\) be the \(q\)-Beckmann boundary from Definition~\ref{def:beckmann-boundary} specialized to $\mathbb N_0^N$ with product Poisson measure $\mu=\operatorname{Pois}(1)^{\otimes N}$, and let \(C_q\) denote the optimal constant determined in Theorems~\ref{thm:bernoulli-l2} and~\ref{thm:exact-phase-transitions}. For \(1\le q\le2\), the Beckmann Poincaré inequality
\[
\|g-\mu(g)\|_{L^1(\mu)}
\le
2^{1/q}C_q\,\mathsf B_q(g),\qquad g\in D(\Delta)
\]
holds on all such spaces.
Moreover, for $q=2$, the constant $\sqrt{\pi}$ is optimal
even for $N=2$.
\end{theorem}

Taking $v=\nabla g$ in the definition of the Beckmann
boundary yields
\[
\|g-\mu(g)\|_{L^1(\mu)}
\le
2^{1/q}C_q\,\mathbb E_\mu|\nabla_qg|.
\]
For $q=2$, this improves the constant \(\sqrt{2\pi}\) obtained from the inequality of Bobkov and Götze to \(\sqrt{\pi}\), while extending the result from Boolean-valued functions to arbitrary functions in \(D(\Delta)\).

\medskip
\noindent\textbf{Finite abelian Cayley graphs.}
Pseudo-Poincar\'e inequalities enter semigroup proofs of Sobolev inequalities; see, for example,~\cite{Ledoux2003}.
Klartag et al.~\cite[Theorem~2.3 and the proof of Lemma~4.1]{KKRT2016} studied pseudo-Poincaré inequalities
on graphs with nonnegative curvature and showed that finite abelian
Cayley graphs satisfy the Bakry--Émery condition \(\mathsf{CD}(0,\infty)\).
In our normalization, their zero-curvature semigroup argument yields
\[
\|g-P_Tg\|_{L^1(\mu)}
\le
\sqrt{2T}\,\mathbb E|\nabla_2g|
\]
for all real-valued functions \(g\).

\begin{theorem}[Abelian Cayley graphs, Corollary~\ref{cor:abelian-cayley-graphs}]
Let \(\mathsf B_2\) be the \(2\)-Beckmann boundary from
Definition~\ref{def:beckmann-boundary} specialized to a finite
abelian Cayley graph \(G=\operatorname{Cay}(\Gamma,S)\) with
the uniform probability measure \(\mu\).
For every \(T>0\), the Beckmann pseudo-Poincaré inequality
\[
\|g-P_Tg\|_{L^1(\mu)}
\le
\sqrt{\frac{4T}{\pi}}\,\mathsf B_2(g),\qquad g:\Gamma\to \mathbb R
\]
holds on all such graphs.
Moreover, the constant
\(\sqrt{4T/\pi}\) is optimal uniformly, even among cyclic Cayley graphs.
\end{theorem}

Taking \(v=\nabla g\) in the definition of the Beckmann boundary yields
\[
\|g-P_Tg\|_{L^1(\mu)}
\le
\sqrt{\frac{4T}{\pi}}\,
\mathbb E|\nabla_2g|.
\]
Thus, our result replaces the classical gradient boundary by the
smaller Beckmann boundary and improves the curvature-based
coefficient from \(\sqrt{2T}\) to \(\sqrt{4T/\pi}\).

\subsection{Proof Strategies}
The proof begins with a semigroup reduction. On a weighted graph, our duality theorem identifies the optimal constant in the Beckmann Poincaré or pseudo-Poincaré inequality with that in a uniform gradient estimate for the integrated heat semigroup. 
In the spirit of the Bismut--Elworthy--Li formula~\cite{Bismut1984, ElworthyLi1994},
we construct processes \((X_t^x)_{t>0}\) and random vector fields \((H_t^x)_{t>0}\) such that
\[
\nabla P_tf(x,y)
=
\mathbb E[f(X_t^x)H_t^x(x,y)],
\]
for all $f\in L^\infty(\mu)$ and $t>0$.
Integrating in time and applying local \(\ell^p\)–\(\ell^q\) duality bounds the optimal constant in the Beckmann inequality by
\[
\sup_{x\in V} \sup_{|v|_q(x)\le 1}\int_0^T\mathbb E\left|\frac{1}{2\mu(x)}\sum_{y\sim x}b(x,y)v(x,y)H_t^x(x,y)\right|\,\mathrm dt.
\]
The problem is thus reduced to bounding the time-integrated directional projections of the Malliavin weights. Crucially, the direction \(v\) remains fixed throughout the time integration.

On the Hamming cube \(\{-1,1\}^N\), the heat process flips the coordinates independently, and the general weight bound is expressed through the integral
\[
\mathfrak B(v)
:=
\int_0^1\mathbb E\left|\sum_{i=1}^N v_i(B_i(s)-s)\right|\frac{\mathrm ds}{s},
\]
where \(B_1(s),\ldots,B_N(s)\) are independent \(\operatorname{Ber}(s)\) random variables.
The cube problem therefore leads us to determine, for \(1\le q\le2\), the dimension-free constants
\[
C_q
:=
\sup_{N\ge1}\sup_{\|v\|_{\ell^q}=1}\mathfrak B(v).
\]

Earlier works also reduced the cube problem to time integrals involving centered Bernoulli variables; see Ivanisvili et al.~\cite{IvanisviliLiVanHandelVolberg2018} and Ivanisvili and Stone~\cite{IvanisviliStone2024}.
In particular, the scalar case of~\cite[Theorem~1.2]{IvanisviliVanHandelVolberg2020} already gives $\|g-\mu(g)\|_{L^1(\mu)}\le \mathbb E\mathfrak B(\nabla g/2)$.
Indeed, the substitution $s=(1-e^{-t})/2$ and the symmetry of centered Bernoulli sums under $s\leftrightarrow1-s$ give the integral defining $\mathfrak B$.
Thus the fixed-coefficient integral reduction is already available before the sign comparison below. The quantitative bounds in~\cite{IvanisviliLiVanHandelVolberg2018,IvanisviliStone2024} estimate the supremum over coefficients inside the time integral.
Applying Cauchy--Schwarz to this integral representation gives the upper bound \(C_2\le\pi/2\), and refinements based on higher moments give a strictly smaller bound.
However, as observed by Ivanisvili and Stone~\cite[Remark~1.3]{IvanisviliStone2024}, even optimal estimates at each fixed time lead to a constant strictly larger than \(\sqrt{\pi/2}\) after integration.

We estimate the fixed-coefficient integral directly to avoid this loss. 
For nonnegative coefficients, the computation of $\mathfrak B(v)$ leads to a binomial-type sum. An identity of Hurwitz~\cite[Eq.~(IV)]{Hurwitz1902}, arising from his study of
Riemann surfaces, yields
\[
\mathfrak B(v)
=
\int_0^\infty\left(1-\prod_{i=1}^N e^{-v_is}(1+v_is)\right)\,\frac{\mathrm ds}{s^2},
\]
which makes the subsequent estimates more tractable. To handle
coefficients of mixed signs, numerical experiments led us to
conjecture the sign-alignment inequality
\[
\mathfrak B(v)\le\mathfrak B(|v|).
\]
The common centering drift also suggests this comparison: coefficients of the same sign reinforce the terms \(-sv_i\), whereas coefficients of opposite signs lead to cancellation.

Our proof uses Fourier analysis twice. The Fourier representation of the absolute-value function first reduces the comparison to a positivity statement involving sine expectations. A second Fourier expansion, now in the time variable \(s\), establishes this positivity mode by mode.
This use of termwise positivity after Fourier expansion also appears in other areas of mathematics. For an example from analytic number theory, see Bondarenko and Seip~\cite[Section~3]{BondarenkoSeip2017}.

To describe the first step, we separate the positive and negative coefficients by setting
\[
X_s^+:=\sum_{v_i>0}v_i(B_i(s)-s),
\qquad
X_s^-:=-\sum_{v_i<0}v_i(B_i(s)-s).
\]
The original sum is \(X_s^+-X_s^-\), while replacing the coefficients by their absolute values gives \(X_s^++X_s^-\). Since $X_s^\pm$ are independent, the characteristic-function formula for the first absolute moment yields
\[
\mathfrak B(|v|)-\mathfrak B(v)
=
\frac4\pi
\int_0^\infty
\left(
\int_0^1
\mathbb E\sin(rX_s^+)\,
\mathbb E\sin(rX_s^-)
\,\frac{\mathrm ds}{s}
\right)
\frac{\mathrm dr}{r^2}.
\]
It is therefore enough to prove that the inner integral is nonnegative for almost every \(r>0\). In fact, we establish this for every \(r>0\).

The second Fourier step relies on a positivity property of the sine kernel. Although the weight \(ds/s\) does not make distinct sine modes orthogonal, it pairs any two nonnegative frequencies nonnegatively
\[
\int_0^1 \sin(as)\sin (bs)\,\frac{\mathrm ds}{s}\ge 0,\qquad a,b\ge 0.
\]
It is thus enough to show that the Fourier sine coefficients of the two expectations have the same sign. For each nonnegative coefficient \(v_i\), the complex Fourier coefficients of the corresponding one-coordinate characteristic function form the nonnegative sequence
\[
p_{v_i}(k)
=
\left(\frac{\sin(rv_i/2)}{\pi k+rv_i/2}\right)^2,
\qquad
k\in\mathbb Z,
\]
with removable singularities understood by continuity.
Independence expresses the coefficients for each Bernoulli sum as convolutions of these sequences. A reflection comparison shows that the convolutions remain biased toward negative frequencies. Hence the associated sine coefficients are all nonpositive, and the positivity of the sine kernel proves the sign-alignment inequality.

The other two models fit into the same framework. For product Poisson spaces, the Malliavin-weight formula again reduces to integrals involving centered Bernoulli sums. For finite abelian Cayley graphs, it instead leads to centered Poisson sums.
We treat the Poisson sums by the same two-step Fourier argument, replacing the Fourier sine series in the Bernoulli case by a Fourier sine integral.

The extremizers also explain the relation between the subquadratic and quadratic results. For $q<2$, the optimal coefficient vectors have finite support and equal nonzero coordinates; the support size increases through the phases and tends to infinity as $q\uparrow2$.
At $q=2$, increasingly many equal coordinates approach the Gaussian constant on the cube. The other sharpness constructions produce Gaussian limits by increasing the initial vertex coordinates in a fixed two-dimensional Poisson space, or the number of generators in the cyclic Cayley graphs. In the Poisson case the measure remains $\operatorname{Pois}(1)^{\otimes2}$; it is the binomial thinning of the initial state that has growing variance.

\subsection{Organization of the paper}
In Section 2, we introduce Beckmann boundaries on weighted graphs, prove a Talagrand-type functional isoperimetric inequality, and establish the Beckmann duality theorem. Section 3 develops gradient estimates for the integrated heat semigroup using Malliavin weights and constructs explicit weights for Hamming cubes, product Poisson spaces, and finite abelian Cayley graphs. In Section 4, we prove Fourier positivity for centered Bernoulli and Poisson sums and deduce the sign-alignment principle. Section 5 establishes sharp \(\ell^q\)-estimates for the Bernoulli functional and a sharp \(\ell^2\)-estimate for the Poisson functional, including the phase transitions in the Bernoulli case. Finally, Section 6 combines these results to derive the Beckmann Poincaré and pseudo-Poincaré inequalities and prove the stated sharpness results.

\section{Beckmann boundaries}

\subsection{Weighted graphs}
A weighted graph is defined as a triple $G=(V,b,\mu)$, where $V$ represents a countable set of vertices, $b:V\times V\to[0,\infty)$ is a symmetric function satisfying $b(x,x)=0$, and $\mu:V\to[0,1]$ is a probability measure on $V$ with full support.
Two vertices $x,y\in V$ are called neighbors, denoted $x\sim y$, if $b(x,y)>0$.
Connectedness means that any two vertices $x,y\in V$ can be joined by a finite path
\[
x=x_0\sim x_1\sim\cdots\sim x_n=y,
\]
while local finiteness requires that every vertex $x\in V$ has finitely many neighbors.
Throughout, we consider only connected and locally finite weighted graphs,
allowing us to define the combinatorial degree and weighted degree of a vertex $x\in V$ as follows
\[
\deg(x):=\#\{y\in V:y\sim x\},
\qquad
\operatorname{Deg}(x):=\frac{1}{\mu(x)}\sum_{y\sim x}b(x,y).
\]

A function on $G$ is a real-valued map $f:V\to\mathbb R$, while a vector field on $G$ is a map $v:V\times V\to\mathbb R$ that vanishes for all non-adjacent pairs of vertices. The associated inner products for these spaces are given by
\[
\langle f,g\rangle_\mu:=\sum_{x\in V}\mu(x)f(x)g(x),\qquad
\langle v,w\rangle_b:=\frac12\sum_{x\in V}\sum_{y\sim x}b(x,y)v(x,y)w(x,y),
\]
and the resulting Hilbert spaces are denoted as $L^2(\mu)$ and $L^2(b/2)$, respectively. Note that no antisymmetry condition is imposed on vector fields.

The edge weights naturally define an energy form on $L^2(\mu)$, which generates the graph Laplacian and its corresponding heat semigroup.
Specifically, we consider the Neumann Dirichlet form $Q$ on $L^2(\mu)$, given by
\[
Q(f,g):=\frac12\sum_{x\in V}\sum_{y\sim x}b(x,y)(f(y)-f(x))(g(y)-g(x))
\]
on its domain $D(Q):=\{f\in L^2(\mu):Q(f,f)<\infty\}$.
The associated non-positive self-adjoint operator \(\Delta\), known as the Neumann realization of the graph Laplacian, acts pointwise as
\[
\Delta f(x)=\frac1{\mu(x)}\sum_{y\sim x}b(x,y)(f(y)-f(x)),
\qquad
f\in D(\Delta),
\qquad 
x\in V.
\]
This operator generates the strongly continuous symmetric Markov semigroup \(P_t:=e^{t\Delta}\), which is conservative in the sense that
\[
P_t1=1,\qquad t\ge 0.
\]
Conservativeness here refers to the Neumann semigroup and should be distinguished from stochastic completeness, which relates to the minimal heat semigroup; see \cite{Haeseler12} for a counterexample.

\begin{assumption}
For every \(f\in L^\infty(\mu)\), we assume the integral
\[
\int_0^\infty(P_tf-\mu(f))\,\mathrm dt
\]
converges in \(L^2(\mu)\), so that the Poisson equation \(\Delta u=f-\mu(f)\) admits a unique mean-zero solution \(u\in D(\Delta)\).
\end{assumption}

To express the Laplacian in divergence form, we introduce the gradient operator \(\nabla:D(Q)\subset L^2(\mu)\to L^2(b/2)\) and define the divergence as its negative Hilbert space adjoint, \(\operatorname{div}:=-\nabla^*\). Their pointwise expressions are
\[
\nabla f(x,y)=f(y)-f(x),
\qquad
\operatorname{div}v(x)=\frac1{2\mu(x)}\sum_{y\sim x}b(x,y)(v(x,y)-v(y,x)).
\]
For the dual argument, we further define the pointwise \(q\)-gradient length of a function \(f\) by $|\nabla_q f|(x):=|\nabla f|_q(x)$,
where
\[
|v|_q(x):=
\begin{cases}
\displaystyle
\left(\frac1{2\mu(x)}
\sum_{y\sim x}b(x,y)|v(x,y)|^q\right)^{1/q}
&1\le q<\infty,\\
\displaystyle\max_{y\sim x}|v(x,y)|,
&q=\infty,
\end{cases}
\]
with empty maximum interpreted as zero. The assumption of local finiteness ensures that both quantities are finite at every vertex.

\subsection{Beckmann boundaries}
The concept of Beckmann boundaries originates from Beckmann's formulation of optimal transport~\cite{Beckmann1952}, where the goal is to minimize the cost of a vector field subject to a prescribed divergence constraint;
see also~\cite[Chapter~4]{Santambrogio2015}.
Motivated by this idea, Ivanisvili, Xie, and Zhang~\cite{IvanisviliXieZhang2026} introduced the notion of the Beckmann boundary on the Hamming cube. We generalize their definition to weighted graphs.

\begin{definition}[$q$-Beckmann boundaries]
\label{def:beckmann-boundary}
For $1\le q\le \infty$ and $f\in D(\Delta)$,
we define the \emph{$q$-Beckmann boundary} of $f$ by
\[
\mathsf B_q(f)
:=
\inf\left\{\mathbb E|v|_q:v\in D(\operatorname{div}),\,\operatorname{div}v=\Delta f\right\}.
\]
\end{definition}

The Beckmann boundary serves as a variational relaxation of the traditional gradient boundary. By choosing \(v=\nabla f\) as an admissible vector field, it is immediate that
\[
\mathsf B_q(f)\le \mathbb E|\nabla_qf|.
\]
The following Talagrand-type inequality gives a perimeter-like lower bound for this boundary. Its variance--supremum-norm formulation differs from the influence-sensitive cube inequalities in \cite{IvanisviliXieZhang2026}.

\begin{theorem}[Talagrand-type functional isoperimetry]\label{thm:talagrand}
Assume that $(P_t)_{t\ge0}$ satisfies the reverse Poincaré inequality
\[
2t|\nabla_2P_th|^2(x)
\le P_t(h^2)(x)-(P_th(x))^2,
\qquad h\in D(\Delta),
\]
and the logarithmic Sobolev inequality
\[
\operatorname{Ent}(h^2)
\le
\frac1\lambda\mathbb E|\nabla_2h|^2,
\qquad
h\in D(Q).
\]
Then, for every $1\le q\le 2$ and every nonconstant $g\in D(\Delta)\cap L^\infty(\mu)$, it holds that
\[
\mathsf B_q(g)
\gtrsim_q
\lambda^{1/q}\frac{\operatorname{Var}(g)}{\|g-\mu(g)\|_{\infty}}
\left(\log\left(1+\frac{\|g-\mu(g)\|_{\infty}^2}{\operatorname{Var}_\mu(g)}\right)\right)^{1/q}.
\]
\end{theorem}

On a finite graph, take $g=\mathbf1_A$ and $a=\min\{\mu(A),1-\mu(A)\}\in(0,1/2]$.
Then $\operatorname{Var}(g)=a(1-a)$ and $\|g-\mu(g)\|_\infty=1-a$, so the theorem gives
\[
\mathsf B_q(\mathbf1_A)\gtrsim_q\lambda^{1/q}a\bigl(\log(1/a)\bigr)^{1/q}.
\]
Thus the variational relaxation retains a perimeter-like isoperimetric scaling.

On finite graphs, the reverse Poincar\'e inequality is equivalent to
the pointwise Bakry--\'Emery curvature-dimension condition
\(\mathsf{CD}(0,\infty)\). Under additional hypotheses, similar equivalences extend to locally finite graphs; see~\cite[Theorem~3]{Hua2019}.
The proof of Theorem~\ref{thm:talagrand} relies on the following lemma.

\begin{lemma}[Heat dissipation]
\label{lem:heat-dissipation}
Assume that $(P_t)_{t\ge0}$ satisfies the reverse Poincar\'e inequality
\[
2t|\nabla_2P_th|^2(x)
\le P_t(h^2)(x)-(P_th(x))^2,
\qquad h\in D(\Delta).
\]
Then, for every $1\le q\le2$ and every nonconstant
$g\in D(\Delta)\cap L^\infty(\mu)$, it holds that
\begin{gather*}
\|g-P_tg\|_{L^1(\mu)}
\lesssim_q t^{1/q}\mathsf B_q(g),\\
\operatorname{Var}_\mu(g)-\operatorname{Var}_\mu(P_tg)
\lesssim_q
t^{1/q}\|g-\mu(g)\|_\infty\mathsf B_q(g).
\end{gather*}
\end{lemma}

\begin{proof}
Let $2\le p\le\infty$ satisfy $1/p+1/q=1$.
For $h\in D(\Delta)\cap L^\infty(\mu)$, the reverse Poincaré inequality combined with the Markov property ensures
\[
\||\nabla_2P_sh|\|_\infty
\le (2s)^{-1/2}\|h\|_\infty,
\qquad
\||\nabla_\infty P_sh|\|_\infty
\le2\|h\|_\infty.
\]
By interpolation, $\||\nabla_pP_sh|\|_\infty\le2^{1-3/p}s^{-1/p}\|h\|_\infty$.

For every $v\in D(\operatorname{div})$ satisfying
$\operatorname{div}v=\Delta g$, self-adjointness and
integration by parts yield
\[
|\langle g-P_tg,h\rangle_\mu|
=
\left|\int_0^t\langle v,\nabla P_sh\rangle_b\,\mathrm ds\right|
\le
2^{1-3/p}t^{1/q}
\left(\int_0^1s^{-1/p}\,\mathrm ds\right)
\|h\|_\infty\mathbb E|v|_q.
\]
The first inequality follows by approximation and duality.
Combining this estimate with self-adjointness, we obtain 
\begin{align*}
\operatorname{Var}_\mu(g)-\operatorname{Var}_\mu(P_tg)
&=
\langle g-P_{2t}g,g-\mu(g)\rangle_\mu\\
&\le
2^{1-3/p}(2t)^{1/q}
\left(\int_0^1s^{-1/p}\,\mathrm ds\right)
\|g-\mu(g)\|_\infty\mathsf B_q(g).
\end{align*}
Rearranging yields the second inequality.
\end{proof}

\begin{proof}[Proof of Theorem~\ref{thm:talagrand}]
Given that $\operatorname{Var}(g)\le \|g-\mu(g)\|_\infty^2$,
we may decompose the logarithmic term
\begin{align*}
\log\left(1+\frac{\|g-\mu(g)\|_\infty^2}{\operatorname{Var}(g)}
\right)
&\le
\log 2+\log\frac{\operatorname{Var}(g)}{\|g-\mu(g)\|_{L^1(\mu)}^2}\\
&\qquad +2\log\frac{\|g-\mu(g)\|_\infty\|g-\mu(g)\|_{L^1(\mu)}}{\operatorname{Var}(g)}.
\end{align*}

To estimate the first and second terms,  we utilize exponential ergodicity to select a time
\[
0<t\le\frac1\lambda\log\sqrt2
\]
such that $\|g-\mu(g)\|_{L^2(\mu)}=\sqrt 2\|P_tg-\mu(g)\|_{L^2(\mu)}$.
For this choice of $t$, hypercontractivity~\cite{Gross1975},
followed by interpolation between $L^1(\mu)$- and $L^2(\mu)$-norms, yields
\begin{align*}
\|g-\mu(g)\|_{L^2(\mu)}
&\le
\sqrt2\|g-\mu(g)\|_{L^{1+e^{-2\lambda t}}(\mu)}\\
&\le
\sqrt2\|g-\mu(g)\|_{L^1(\mu)}^{\tanh(\lambda t)}\|g-\mu(g)\|_{L^2(\mu)}^{1-\tanh(\lambda t)}.
\end{align*}
Rearranging terms, taking logarithms, and applying Lemma~\ref{lem:heat-dissipation}, we obtain
\begin{align*}
\log 2+\log\frac{\operatorname{Var}(g)}
{\|g-\mu(g)\|_1^2}
&\le
\log 2+\frac{\log2}{\tanh(\lambda t)}\lesssim\frac1{\lambda t}
\lesssim_q\left(\frac{\|g-\mu(g)\|_\infty\,\mathsf B_q(g)}{\lambda^{1/q}\operatorname{Var}(g)}\right)^q.
\end{align*}

For the third term, we again use exponential ergodicity to select a time 
\[
0<s\le\frac1\lambda\left(\log2+\frac12
\log\frac{\operatorname{Var}(g)}{\|g-\mu(g)\|_{L^1(\mu)}^2}\right)
\]
such that $\|g-\mu(g)\|_{L^1(\mu)}=2\|P_sg-\mu(g)\|_{L^1(\mu)}$.
For this choice of $s$,
Lemma~\ref{lem:heat-dissipation} gives
\begin{align*}
\log\frac{\|g-\mu(g)\|_\infty\|g-\mu(g)\|_1}
{\operatorname{Var}(g)}
&\lesssim_q
\left(\frac{\|g-\mu(g)\|_\infty\|g-\mu(g)\|_1}{\operatorname{Var}(g)}\right)^{q/2}\\
&\lesssim_q
s^{1/2}\left(\frac{\|g-\mu(g)\|_\infty\mathsf B_q(g)}{\operatorname{Var}(g)}\right)^{q/2}\\
&\lesssim_q
\left(\frac{\|g-\mu(g)\|_\infty\,\mathsf B_q(g)}{\lambda^{1/q}\operatorname{Var}(g)}\right)^q.
\end{align*}This completes the proof.
\end{proof}

\subsection{Duality}
The Beckmann Poincar\'e and pseudo-Poincar\'e inequalities admit equivalent dual formulations, expressed as uniform gradient estimates for the integrated heat semigroup.
We use the convention $P_\infty g:=\mu(g)$.

\begin{theorem}[Beckmann duality]
\label{thm:beckmann-duality}
Let $0<T\le \infty$ be fixed, let $1\le p,q\le \infty$ satisfy $1/p+1/q=1$, and define the integrated heat semigroup \(R_T\) by
\[
R_Tf:=\int_0^T (P_sf-\mu(f))\,\mathrm ds\qquad \text{in $L^2(\mu)$}.
\]
For every $C>0$, the following assertions are equivalent:
\begin{enumerate}
\item For every $f\in L^\infty(\mu)$, it holds that
\[
\||\nabla_{p}R_Tf|\|_{\infty}
\le
C\|f\|_{\infty}.
\]

\item For every $g\in D(\Delta)$, it holds that 
\[
\|g-P_Tg\|_{L^1(\mu)}
\le
C\mathsf B_{q}(g).
\]
\end{enumerate}

\end{theorem}

\begin{proof}
For bounded $f$, $R_Tf\in D(\Delta)$ and $\Delta R_Tf=P_Tf-f$, also when $T=\infty$ by the standing assumption and closedness of $\Delta$.
We start by recording two key duality identities that facilitate the connection between the two conditions:
\begin{gather*}
\|g-P_Tg\|_{L^1(\mu)}
=
\sup\{|\langle g-P_Tg,f\rangle_\mu|:\text{$f\in L^\infty(\mu)$, $\|f\|_{\infty}\le 1$}\},\\
\||\nabla_pR_Tf|\|_{\infty}
=
\sup\{|\langle v,\nabla R_Tf\rangle_b|:\text{$v$ finitely supported, $\mathbb E|v|_q\le 1$}\}.
\end{gather*}
The first identity follows from the standard duality of $L^1(\mu)$ and $L^\infty(\mu)$, while the second stems from the duality between $L^1(\mu;\ell^q)$ and $L^\infty(\mu;\ell^p)$.

\textbf{Step 1.}
Assume that condition (1) holds.
Let $g\in D(\Delta)$ be fixed, and choose a vector field $v\in D(\operatorname{div})$ such that $\operatorname{div}v=\Delta g$.
For every $f\in L^\infty(\mu)$ with $\|f\|_{\infty}\le1$,
self-adjointness of \(P_T\) and \(\Delta\), together with the relation $\operatorname{div}=-\nabla^*$, leads us to
\[
\langle g-P_Tg,f\rangle_\mu
=
-\langle\operatorname{div}v,R_Tf\rangle_\mu
=
\langle v,\nabla R_Tf\rangle_b.
\]
Applying the uniform gradient estimate for the integrated heat semigroup in condition (1), we can bound the right-hand side
\[
|\langle g-P_Tg,f\rangle_\mu|
\le
\mathbb E\left[|v|_q|\nabla_pR_Tf|\right]\le
C\mathbb E|v|_q.
\]
Taking the supremum over all \(f\in L^\infty(\mu)\) with $\|f\|_\infty\le 1$,
and then the infimum over all \(v\in D(\operatorname{div})\) with \(\operatorname{div}v=\Delta g\), proves condition (2).

\textbf{Step 2.}
Now assume that condition (2) holds.
Let $f\in L^\infty(\mu)$ be fixed, and choose a finitely supported vector field $v$ such that $\mathbb E|v|_q\le 1$.
Then $v\in D(\operatorname{div})$, and $\operatorname{div}v$ has finite support, is bounded, and has mean zero.
The standing assumption therefore provides a mean-zero $g\in D(\Delta)$ with $\Delta g=\operatorname{div}v$, for which a similar argument gives
\[
\langle g-P_Tg,f\rangle_\mu
=
-\langle\operatorname{div}v,R_Tf\rangle_\mu
=
\langle v,\nabla R_Tf\rangle_b.
\]
Using condition (2), we can estimate
\[
|\langle v,\nabla R_Tf\rangle_b|
\le
C\|f\|_{L^\infty(\mu)}\mathsf B_q(g)
\le
C\|f\|_{L^\infty(\mu)}\mathbb E|v|_q.
\]
Taking the supremum over all finitely supported vector fields \(v\) with \(\mathbb E_\mu|v|_q\le1\) establishes condition (1).
\end{proof}

Theorem~\ref{thm:beckmann-duality} reduces the problem of determining the optimal constants in Beckmann Poincar\'e and pseudo-Poincar\'e inequalities to establishing uniform gradient estimates for the integrated heat semigroup \(R_T\), without compromising the constants. The variational definition of \(\mathsf B_q(g)\) further guarantees the smallest possible boundary term within this duality framework, achieved by optimizing over all vector fields \(v\) satisfying \(\operatorname{div}v=\Delta g\). Since \(\mathsf B_q(g)\le\mathbb E_\mu|\nabla_q g|\), this result also recovers the classical inequalities with the same constant. However, pinpointing the optimal constants in these classical inequalities remains a subtle problem.

\section{Malliavin weights}
\label{sec:malliavin-weights}

\subsection{Gradient estimates}

The following proposition mirrors the classical
Bismut--Elworthy--Li formula~\cite{Bismut1984,ElworthyLi1994},
which expresses spatial derivatives of the heat semigroup as expectations involving Malliavin weights.

\begin{proposition}[Malliavin weights]
\label{prop:weight-time-dependent}
Fix $x\in V$, and let $(X_t^x)_{t>0}$ be a jointly measurable \(V\)-valued process satisfying
\[
\mathbb E[f(X_t^x)]=P_tf(x)
\]
for all $f\in L^\infty(\mu)$ and $t>0$.
Then there exists a jointly measurable vector-field-valued process \((H_t^x)_{t>0}\) on \(G\) such that
\[
\nabla P_tf(x,\cdot)=\mathbb E[f(X_t^x)H_t^x(x,\cdot)]
\]
for all $f\in L^\infty(\mu)$ and $t>0$.
\end{proposition}

\begin{proof}
Let \(p_t(y,z)\) be the heat kernel of \(P_t\) with respect to \(\mu\), so that
\[
P_tf(y)=\sum_{z\in V}p_t(y,z)f(z)\mu(z).
\]
Since the graph is connected, the kernel $p_t(y,z)$ is strictly positive for all $t>0$ and $y,z\in V$; see \cite[Theorem~7.3]{Haeseler12}. Define the random vector field \(H_t\) as:
\[
H_t(y,z)
:=
\frac{p_t(z,X_t^x)-p_t(y,X_t^x)}
{p_t(x,X_t^x)},
\qquad \text{for }y\sim z,
\]
and set $H_t(y,z)=0$ otherwise.
Joint measurability follows from that of $(X^x_t)_{t>0}$ and the heat kernel $(p_t)_{t>0}$. Using the distribution of \(X_t^x\), we compute
\[
\mathbb E[f(X_t^x)H_t(x,y)]
=
\sum_{z\in V}f(z)(p_t(y,z)-p_t(x,z))\mu(z)
=
\nabla P_tf(x,y),
\]
which holds for all $f\in L^\infty(\mu)$ and $t>0$.
\end{proof}

This next proposition relates gradient estimates for the integrated heat semigroup to time-integrated bounds on directional projections of the Malliavin weights.

\begin{proposition}[Gradient estimates for the integrated heat semigroup]
\label{prop:heat-gradient-estimate}
Let $0<T\le\infty$, and let $1\le p,q\le\infty$ satisfy $1/p+1/q=1$.
For the Malliavin weights $(H_t)_{t>0}$ defined in Proposition~\ref{prop:weight-time-dependent}, we have
\[
\sup_{\|f\|_\infty\le1}\||\nabla_pR_Tf|\|_\infty
\le
\sup_{x\in V}\sup_{|v|_q(x)\le1}\int_0^T\mathbb E\left|\frac1{2\mu(x)}\sum_{y\sim x}b(x,y)v(x,y)H_t(x,y)\right|\,\mathrm dt.
\]
\end{proposition}

\begin{proof}
By the definition of the integrated heat semigroup $R_T$,
\[
\nabla R_Tf(x,y)
=
\int_0^T\nabla P_tf(x,y)\,\mathrm dt
=
\int_0^T\mathbb E[f(X_t^x)H_t(x,y)]\,\mathrm dt.
\]
For finite truncations this follows by linearity. When $T=\infty$, the integrals are understood as improper limits: the standing assumption gives $L^2(\mu)$ convergence, hence convergence at every vertex since $\mu(x)>0$, and therefore along each edge.
Hence, by local $\ell^p$--$\ell^q$ duality,
\begin{align*}
\||\nabla_pR_Tf|\|_\infty
&=\sup_{x\in V}
\sup_{|v|_q(x)\le1}
\left|
\frac1{2\mu(x)}
\sum_{y\sim x}b(x,y)v(x,y)\nabla R_Tf(x,y)
\right|\\
&\le
\|f\|_\infty\sup_{x\in V}
\sup_{|v|_q(x)\le1}
\int_0^T
\mathbb E\left|
\frac1{2\mu(x)}
\sum_{y\sim x}b(x,y)v(x,y)H_s(x,y)
\right|\,\mathrm ds.
\end{align*}
Taking the supremum over $\|f\|_\infty\le1$ proves the assertion.
\end{proof}

\subsection{Explicit weights}

Proposition~\ref{prop:weight-time-dependent} demonstrates the existence
of Malliavin weights representing $\nabla P_tf(x,\cdot)$.
Below, we present explicit constructions of such weights for specific probabilistic structures of the semigroup in three illustrative examples.

\begin{example}[Hamming cubes]
\label{ex:cube-weights}
Let $G=\{-1,1\}^N$ represent the $N$-dimensional Hamming cube, equipped with the uniform probability measure $\mu$. The Laplacian is defined as
\[
\Delta f(x)=\frac12\sum_{y\sim x}(f(y)-f(x)).
\]

For $t>0$, let $B_1(s),\dots,B_N(s)$ be independent Bernoulli random variables satisfying
\[
B_i(s)\sim \operatorname{Ber}(s),
\qquad
s:=\frac{1-e^{-t}}{2}.
\]
By~\cite[Lemma~2.1]{IvanisviliVanHandelVolberg2020},
we can construct the process $(X_t)_{t>0}$ and the Malliavin weights $(H_t)_{t>0}$ as follows
\[
X_t^x:=(x_i(-1)^{B_i(s)})_{i=1}^N,
\qquad
H_t(x)=\frac{e^{-t}}{s(1-s)}(B_i(s)-s)_{i=1}^N.
\]
Indeed, flipping the initial coordinate replaces $B\sim\operatorname{Ber}(s)$ by $1-B$. For $b\in\{0,1\}$,
\[
\frac{\mathbb P(1-B=b)}{\mathbb P(B=b)}-1
=\frac{1-2s}{s(1-s)}(b-s)=\frac{e^{-t}}{s(1-s)}(b-s).
\]
This likelihood-ratio weight represents the difference $P_tf(x^{\oplus i})-P_tf(x)$ and explains the centered Bernoulli term.
Consequently, the optimal constant in the associated Beckmann inequality can be bounded above by
\[
\frac12
\sup_{x\in V}\sup_{|v|_q(x)\le1}
\left(\int_0^{(1-e^{-T})/2}+\int_{(1+e^{-T})/2}^1\right)
\mathbb E\left|\sum_{i=1}^N v_i(B_i(s)-s)\right|\,\frac{\mathrm ds}{s}.
\]
Here \(|v|_q(x)\) denotes the local weighted norm, which does not coincide with the standard \(\ell^q\)-norm of the coefficient vector.
\end{example}

\begin{example}[Product Poisson spaces]
\label{ex:poisson-weights}
Let $G=\mathbb N_0^N$ be equipped with the product Poisson measure
$\mu=\operatorname{Pois}(1)^{\otimes N}$.
Join two vertices whenever they differ by $\pm e_i$ for some $1\le i\le N$,
where $e_i$ denotes the $i$-th standard basis vector of $\mathbb R^N$.
The Laplacian on $G$ is given by
\[
\Delta f(x)
=
\sum_{i=1}^N(f(x+e_i)-f(x))+\sum_{x_i\ge1}x_i(f(x-e_i)-f(x)).
\]

For $t>0$ and $x\in V$, let $P_1(s),\dots,P_N(s)$ and $B_1^x(s),\dots,B_N^x(s)$ be independent random variables such that
\[
P_i(s)\sim\operatorname{Pois}(s),
\qquad
B_i^x(s)\sim\operatorname{Bin}(x_i,s),
\qquad 
s:=1-e^{-t}.
\]
Using these, we may choose the process $(X_t)_{t>0}$ and the Malliavin weights $(H_t)_{t>0}$ as follows
\begin{gather*}
X_t^x=x+\left(P_i(s)-B_i^x(s)\right)_{i=1}^N,\\
H_t(x,x-e_i)=\frac{B_i^x(s)-x_is}{x_is},\qquad H_t(x,x+ e_i)=\frac{(1-s)(P_i(s)-s)}{s}.
\end{gather*}
The backward formula is used only when $x_i\ge1$.
Indeed, by applying the Poisson and binomial integration-by-parts formulas, we confirm that
\begin{gather*}
\nabla P_tf(x,x-e_i)=\frac1{x_is}\mathbb E[B_i^x(s)f(X_t^x)]-\mathbb E[f(X_t^x)]=\mathbb E\left[
f(X_t^x)\frac{B_i^x(s)-x_is}{x_is}
\right],\\
\nabla P_tf(x,x+e_i)
=(1-s)\mathbb E[f(X_t^x+e_i)-f(X_t^x)]=\mathbb E\left[f(X_t^x)\frac{(1-s)(P_i(s)-s)}{s}\right].
\end{gather*}
As a result, the optimal constant in the corresponding Beckmann inequality can be bounded above by
\begin{align*}
\frac12\sup_{x\in V}\sup_{|v|_q(x)\le1}\int_0^{1-e^{-T}}\mathbb E\,&\Bigg|(1-s)\sum_{i=1}^N v(x,x+e_i)(P_i(s)-s)\\
&\qquad +\sum_{x_i\ge1}v(x,x-e_i)(B_i^x(s)-x_is)\Bigg|\,\frac{\mathrm ds}{s(1-s)}.
\end{align*}

Note that this bound reduces to an estimate for centered Bernoulli sums. Since the preceding integral depends only on the joint law at each fixed \(s\), we may realize the \(P_i\) as unit-rate Poisson processes, preserving the independence of the full family. By independent increments,
\[
\mathbb E[P_i(s)-sP_i(1)\mid (P_j(s),B_j^x(s))_{j=1}^N]
=(1-s)(P_i(s)-s).
\]
Thus, by conditional Jensen’s inequality and the tower property,
the expectation in the preceding integral is bounded above by
\[
\mathbb E\left|\sum_{i=1}^N v(x,x+e_i)(P_i(s)-sP_i(1)) +\sum_{x_i\ge1}v(x,x-e_i)(B_i^x(s)-x_i s)\right|,
\]
and the optimal constant is bounded above by
\begin{align*}
\frac12\sup_{x\in V}\sup_{|v|_q(x)\le1}
\left(\int_0^{1-e^{-T}}+\int_{e^{-T}}^1\right)
\mathbb E\Bigg|
&\sum_{i=1}^N v(x,x+e_i)(P_i(s)-sP_i(1))\\
&\qquad +\sum_{x_i\ge1}v(x,x-e_i)(B_i^x(s)-x_is)\Bigg|\,\frac{\mathrm ds}{s}.
\end{align*}

\end{example}

\begin{example}[Abelian Cayley graphs]
\label{ex:abelian-cayley-weights}
Let \(G=\operatorname{Cay}(\Gamma,S)\) be a finite abelian Cayley graph, where \(\Gamma\) is a finite abelian group and \(S=\{a_1,\ldots,a_N\}\subset\Gamma\setminus\{0\}\) is a symmetric generating set.
Equip $G$ with the uniform probability measure \(\mu\), and join two vertices by an edge whenever they differ by an element of $S$.
The Laplacian on $G$ is given by
\[
\Delta f(x)
=
\frac12\sum_{i=1}^N(f(x+a_i)-f(x)).
\]

For $t>0$, let $P_1(s),\dots,P_N(s)$ be independent
Poisson random variables with parameter $s$, where
\[
P_i(s)\sim\operatorname{Pois}(s),
\qquad s:=\frac t2.
\]
Using this, we can explicitly define the process $(X_t)_{t>0}$ and the Malliavin weights $(H_t)_{t>0}$ as 
\[
X_t^x:=x+\sum_{i=1}^N P_i(s)a_i,
\qquad
H_t(x,x+a_i):=\frac{P_i(s)-s}{s}.
\]
Indeed, the validity of these Malliavin weights can be verified using the Poisson integration-by-parts identity
\[
\nabla P_tf(x,x+a_i)=\mathbb E[f(X_t^x+a_i)-f(X_t^x)]=\mathbb E\left[
f(X_t^x)\frac{P_i(s)-s}{s}
\right].
\]
Therefore, the optimal constant in the associated Beckmann inequality can be bounded above by
\[
\frac12\sup_{x\in V}\sup_{|v|_q(x)\le1}\int_0^{T/2}\mathbb E\left|\sum_{i=1}^N v_i(P_i(s)-s)\right|\,\frac{\mathrm ds}s.
\]
\end{example}

The preceding examples bound the corresponding optimal constants by integrals involving centered Bernoulli and Poisson sums. In the next two sections, we develop the integral estimates for centered Bernoulli and Poisson sums needed to bound the optimal constants in the preceding examples.

\section{Fourier positivity}
\label{sec:fourier-positivity}

The integrals in Section~\ref{sec:malliavin-weights} lead to the introduction of the following functionals.
Consider independent random variables $B_1(s),\dots,B_N(s)$ (resp. $P_1(s),\dots,P_N(s)$) that follow the $\operatorname{Ber}(s)$ distribution (resp. the $\operatorname{Pois}(s)$ distribution).
For $v\in\mathbb R^N$ and $T>0$, define
\begin{align*}
\mathfrak B(v)
:=\int_0^1
\mathbb E\left|\sum_{i=1}^N v_i(B_i(s)-s)\right|
\,\frac{\mathrm ds}{s},\qquad 
\mathfrak P_T(v)
:=\int_0^T
\mathbb E\left|\sum_{i=1}^N v_i(P_i(s)-s)\right|
\,\frac{\mathrm ds}{s}.
\end{align*}
Our goal is to show that replacing the coefficients by their
absolute values can only increase these functionals:
\[
\mathfrak B(v)\le\mathfrak B(|v|),
\qquad
\mathfrak P_T(v)\le\mathfrak P_T(|v|),
\]
where $|v|=(|v_1|,\dots,|v_N|)$.
We establish these comparisons through a Fourier argument, starting with the following elementary lemma.

\begin{lemma}
\label{lem:sine-kernel}
For every $a,b\ge 0$ and $T>0$, it holds that
\[
\int_0^T\sin(as)\sin(bs)\,\frac{\mathrm ds}s\ge0.
\]
Moreover, equality holds if and only if $ab=0$.
\end{lemma}

\begin{proof}
By the product-to-sum formula and a change of variables,
\[
\int_0^T\sin(as)\sin(bs)\,\frac{\mathrm ds}s
=
\frac12\int_{|a-b|T}^{|a+b|T}\frac{1-\cos r}{r}\,\mathrm dr
\ge0.
\]
The last integral vanishes precisely when $|a+b|=|a-b|$,
or equivalently, $ab=0$.
\end{proof}

\subsection{Bernoulli positivity}
We establish Fourier positivity for centered Bernoulli sums.
To prove it, we expand their sine expectations in Fourier sine series and show that the coefficients are nonpositive.

\begin{theorem}[Bernoulli positivity]
\label{thm:bernoulli-positivity}
Let $v\in\mathbb R^N, w\in\mathbb R^M$ with $v_i,w_j\ge 0$. For $0<s<1$, let
\[
B_1(s),\dots,B_N(s),C_1(s),\dots,C_M(s)\sim\operatorname{Ber}(s)
\]
be independent Bernoulli random variables.
Define $X_s:=\sum_{i=1}^N v_i(B_i(s)-s)$ and $Y_s:=\sum_{j=1}^M w_j(C_j(s)-s)$.
Then, for every $r>0$,
\[
\int_0^1\mathbb E\sin(rX_s)\mathbb E\sin(rY_s)\,\frac{\mathrm ds}{s}
\ge
0.
\]
Moreover, equality holds if and only if at least one of the two random sums vanishes identically.
\end{theorem}

Lemma~\ref{lem:sine-kernel} shows that the basic sine modes have nonnegative pairings with respect to the measure \(\mathrm ds/s\).
One of our ideas is therefore to expand \(s\mapsto\mathbb E\sin(rX_s)\) and \(s\mapsto\mathbb E\sin(rY_s)\) in Fourier sine series 
\[
\mathbb E\sin(rX_s)
=
\sum_{n=1}^\infty c_X(n)\sin(2\pi ns),
\qquad
\mathbb E\sin(rY_s)
=
\sum_{m=1}^\infty c_Y(m)\sin(2\pi ms),
\]
and compare their Fourier coefficients.

In practice,
we expand \(s\mapsto\mathbb Ee^{irX_s}\) and \(s\mapsto\mathbb Ee^{irY_s}\) in complex Fourier series and take their imaginary parts. Since $X_s$ and $Y_s$ have the same form, we describe the construction only for \(X_s\).
For $N=1$, a direct calculation gives
\[
\mathbb Ee^{irv(B(s)-s)}
=
\sum_{k\in\mathbb Z}p_v(k)e^{i2\pi ks},
\qquad
p_v(k)
=
\left(\frac{\sin(rv/2)}{\pi k+rv/2}\right)^2.
\]
Removable singularities are interpreted by continuity.
For general \(N\ge 1\), independence therefore gives
\[
\mathbb E e^{irX_s}
=
\sum_{k\in\mathbb Z}p_X(k)e^{i2\pi ks},
\qquad
p_X
:=
p_{v_1}*\cdots*p_{v_N}.
\]
It follows that $c_X(n)=p_X(n)-p_X(-n)$.

\begin{lemma}[Negative biasedness]
\label{lem:reflection-inequality}
For every $n\ge 1$,
\[
p_X(n)\le p_X(-n).
\]
\end{lemma}

\begin{proof}
For each \(1\le i\le N\), decompose $p_{v_i}$ into its symmetric part and the remaining one-sided mass by setting
\[
s_{v_i}(n):=p_{v_i}(|n|), \qquad r_{v_i}(n):=p_{v_i}(n)-p_{v_i}(|n|).
\]
Then \(s_{v_i}\) is symmetric unimodal, \(r_{v_i}\) is nonnegative with support on the negative integers, and 
\[
p_X=(s_1+r_1)*\cdots*(s_N+r_N).
\]
After expansion, each term has the form \(S*R\), with \(S\) symmetric unimodal and \(R\) nonnegative and supported on the nonpositive integers.
For $n\ge1$, every such convolution \(S*R\) satisfies
\[
(S*R)(n)-(S*R)(-n)=\sum_{k\in\mathbb Z}(S(n-k)-S(-n-k))R(k)\le 0,
\]
and hence, summing over all terms in the expansion, $p_X(n)\le p_X(-n)$.
\end{proof}

With the negative biasedness property in hand, we prove the Fourier positivity theorem.

\begin{proof}[Proof of Theorem \ref{thm:bernoulli-positivity}]

Each $p_{v_i}(k)=O((1+|k|)^{-2})$, and this bound is preserved under finite convolution: split the sum according as $|k|\ge|n|/2$ or $|n-k|\ge|n|/2$ and use the $\ell^1$ norm of the other factor. Hence $c_X(n),c_Y(n)=O(n^{-2})$.
Together with Lemma~\ref{lem:reflection-inequality}, this gives $c_X(n),c_Y(m)\le0$ and
\[
\sum_{n=1}^\infty\sqrt n\,|c_X(n)|<\infty,
\qquad
\sum_{m=1}^\infty\sqrt m\,|c_Y(m)|<\infty.
\]
Using $|\sin(2\pi ns)|\le\sqrt{2\pi ns}$, we may therefore apply Fubini's theorem to obtain 
\[
\int_0^1
\mathbb E\sin(rX_s)\mathbb E\sin(rY_s)
\,\frac{\mathrm ds}{s}=
\sum_{n,m\ge1}
c_X(n)c_Y(m)
\int_0^1
\sin(2\pi ns)\sin(2\pi ms)\,\frac{\mathrm ds}s
\ge
0.
\]
This proves the inequality. 

For the equality case, suppose that both $v$ and $w$ have at least one nonzero coordinate.
Differentiating at $s=0$ gives
\begin{gather*}
\left.\frac{\mathrm d}{\mathrm ds}\mathbb E\sin(rX_s)\right|_{s=0}
=
\sum_{i=1}^N(\sin(rv_i)-rv_i)<0,\\
\left.\frac{\mathrm d}{\mathrm ds}
\mathbb E\sin(rY_s)\right|_{s=0}
=
\sum_{j=1}^M(\sin(rw_j)-rw_j)<0,
\end{gather*}
so there exist integers $n_0,m_0\ge1$ such that $c_X(n_0),c_Y(m_0)<0$.
It follows that
\[
\int_0^1 \mathbb E\sin(rX_s)\mathbb E\sin(rY_s)\,\frac{\mathrm ds}{s} 
\ge 
c_X(n_0)c_Y(m_0)\int_0^1 \sin(2\pi n_0s)\sin(2\pi m_0s)\,\frac{\mathrm ds}{s}
>
0.
\]
This proves the rigidity statement.
\end{proof}

\subsection{Poisson positivity}
We next establish the analogous positivity result for centered
Poisson sums with nonnegative coefficients.
Here we use Fourier sine integrals and show that the coefficient
functions are nonpositive.

\begin{theorem}[Poisson positivity]
\label{thm:poisson-positivity}
Let $v\in\mathbb R^N, w\in\mathbb R^M$ with $v_i,w_j\ge 0$. For $s>0$, let
\[
P_1(s),\dots,P_N(s),Q_1(s),\dots,Q_M(s)\sim \operatorname{Pois}(s)
\] be independent Poisson random variables.
Define $X_s:=\sum_{i=1}^N v_i(P_i(s)-s)$ and 
$Y_s:=\sum_{j=1}^M w_j(Q_j(s)-s)$.
Then, for every $r>0$ and $0<T<\infty$, 
\[
\int_0^T \mathbb E\sin(rX_s)\mathbb E\sin(rY_s)
\,\frac{\mathrm ds}s
\ge
0.
\]
Moreover, equality holds if and only if at least one of the two random sums vanishes identically.
\end{theorem}

As in the Bernoulli case, one of our ideas is therefore to expand \(s\mapsto\mathbb E\sin(rX_s)\) and \(s\mapsto\mathbb E\sin(rY_s)\) in sine modes and compare their coefficients.
For the Poisson model, the natural expansions are Fourier sine integrals
\[
\mathbb E\sin(rX_s)=\int_0^\infty c_X(\xi)\sin(\xi s)\,\mathrm d\xi,
\qquad
\mathbb E\sin(rY_s)=\int_0^\infty c_Y(\eta)\sin(\eta s)\,\mathrm d\eta.\]
To identify these coefficients, set
\begin{gather*}
a_X:=\sum_{i=1}^N(1-\cos(rv_i)),
\qquad
b_X:=\sum_{i=1}^N(rv_i-\sin(rv_i)),\\
a_Y:=\sum_{j=1}^M(1-\cos(rw_j)),
\qquad
b_Y:=\sum_{j=1}^M(rw_j-\sin(rw_j)).
\end{gather*}
All four quantities are nonnegative.
When \(a_X,a_Y>0\), a direct computation of the Fourier sine coefficients gives
\begin{gather*}
c_X(\xi)
=
-\frac{4a_Xb_X\xi}
{\pi(a_X^2+(b_X-\xi)^2)
(a_X^2+(b_X+\xi)^2)}
\le 0,\\
c_Y(\eta)
=
-\frac{4a_Yb_Y\eta}
{\pi(a_Y^2+(b_Y-\eta)^2)
(a_Y^2+(b_Y+\eta)^2)}
\le 0.
\end{gather*}

\begin{proof}[Proof of Theorem~\ref{thm:poisson-positivity}]
When \(a_X,a_Y>0\), the formulas above give \(c_X,c_Y\le0\) and
\[
\int_0^\infty \xi|c_X(\xi)|\,\mathrm d\xi<\infty,
\qquad
\int_0^\infty \eta|c_Y(\eta)|\,\mathrm d\eta<\infty.
\]
Using \(|\sin u|\le|u|\), we may therefore apply Fubini's theorem to obtain 
\begin{align*}
&\int_0^T\mathbb E\sin(rX_s)\mathbb E\sin(rY_s)\,\frac{\mathrm ds}{s}\\
&\qquad=
\int_0^\infty\int_0^\infty c_X(\xi)c_Y(\eta)
\left(\int_0^T\sin(\xi s)\sin(\eta s)\,\frac{\mathrm ds}{s}\right)\,\mathrm d\xi\,\mathrm d\eta
\ge0.
\end{align*}
When \(a_X=0\) or \(a_Y=0\), the same argument applies with the corresponding Fourier sine integral replaced by a single sine mode.

The rigidity statement follows as in the proof of Theorem~\ref{thm:bernoulli-positivity}, with strict negativity of a single coefficient replaced by strict negativity of the coefficient function on a set of positive measure.
\end{proof}

\subsection{Sign alignment}
For every \(v\in\mathbb R^N\), we write
\[
v^+:=(v_i^+)_{i=1}^N,
\qquad
v^-:=(v_i^-)_{i=1}^N,
\qquad
|v|:=(|v_i|)_{i=1}^N,
\]
where $v_i^+:=v_i\vee0$ and $v_i^-:=(-v_i)\vee0$.
The Fourier positivity results established above allow us to compare
the values of $\mathfrak B$ and $\mathfrak P_T$ at $v$ and at $|v|$.

\begin{theorem}[Sign alignment]
\label{thm:sign-alignment}
For every $v\in\mathbb R^N$ and $T>0$, it holds that
\[
\mathfrak B(v)\le\mathfrak B(|v|),
\qquad
\mathfrak P_T(v)\le\mathfrak P_T(|v|).
\]
Moreover, equality holds in either inequality if and only if \(v^+=0\) or \(v^-=0\), equivalently, all nonzero coordinates of \(v\) have the same sign.
\end{theorem}

\begin{proof}
We prove the Bernoulli case.
Splitting the coefficients according to their sign, we write
\[
X_s^+:=\sum_{i=1}^N v_i^+(B_i(s)-s),
\qquad
X_s^-:=\sum_{i=1}^N v_i^-(B_i(s)-s).
\]
Since these two sums are independent,
the characteristic-function formula for the first absolute moment gives
\[
\mathfrak B(|v|)-\mathfrak B(v)
=
\frac4\pi\int_0^\infty\left(\int_0^1\mathbb E\sin(rX_s^+)\mathbb E\sin(rX_s^-)\,\frac{\mathrm ds}{s}\right)\,\frac{\mathrm dr}{r^2}
\ge
0.
\]
Equality is immediate if \(v^+=0\) or \(v^-=0\).
Otherwise, the inner integral
is strictly positive for every \(r>0\), and hence
\(\mathfrak B(|v|)-\mathfrak B(v)>0\).
\end{proof}

\section{Sharp Bernoulli and Poisson inequalities}
\label{sec:sharp-inequalities}

In this section, we derive the optimal constants for the $\ell^q$ estimates of $\mathfrak B$ for $1\le q\le2$, as well as the $\ell^2$ estimate of $\mathfrak P_T$. By leveraging the results of Theorem~\ref{thm:sign-alignment}, it is sufficient to consider nonnegative coefficients.

\subsection{Sharp Bernoulli inequalities}
We derive an explicit integral representation of $\mathfrak B$ for
nonnegative coefficients, which yields the sharp $\ell^2$-estimate.

\begin{proposition}
Let $v\in\mathbb R^N$ with $v_i\ge0$. Then
\[
\mathfrak B(v)
=
\int_0^\infty \left(1-\prod_{i=1}^N e^{-v_i s}(1+v_i s)\right)\,\frac{\mathrm ds}{s^2}.
\]
\end{proposition}

\begin{proof}
The case $v=0$ is immediate.
By homogeneity, we may assume that $\sum_{i=1}^Nv_i=1$, and set 
\[
X_t:=\sum_{i=1}^Nv_i(B_i(t)-t)=\sum_{i=1}^N v_iB_i(t)-t.
\]
By the symmetry $X_{1-t}\overset d=-X_t$, the identity
$\mathbb EX_t=0$, and Bernoulli integration by parts, we obtain
\begin{align*}
\mathfrak B(v)
&=
\sum_{i=1}^N v_i\int_0^1\mathbb E[(B_i(t)-t)\mathbf1_{\{X_t>0\}}]\,\frac{\mathrm dt}{t(1-t)}\\
&=
\sum_{i=1}^N v_i\int_0^1
\mathbb E\left[
\mathbf1_{\{\sum_{j\ne i}v_jB_j(t)+v_i>t\}}
-
\mathbf1_{\{\sum_{j\ne i}v_jB_j(t)>t\}}
\right]\,\mathrm dt\\
&=
\sum_{i=1}^N v_i\int_0^1
\mathbb P\left(
\sum_{j\ne i}v_jB_j(t)\le t
<
\sum_{j\ne i}v_jB_j(t)+v_i
\right)\,\mathrm dt,
\end{align*}
where the last equality uses $v_i\ge0$.
Put $J_i=\{1,\ldots,N\}\setminus\{i\}$ for $1\le i\le N$
and $v_A=\sum_{j\in A}v_j$ for $A\subseteq J_i$.
Expanding the probabilities and translating each integration interval, we obtain
\begin{align*}
\mathfrak B(v)
&=
\sum_{i=1}^N v_i\int_0^1
\mathbb P\left(
\sum_{j\ne i}v_jB_j(t)\le t<
\sum_{j\ne i}v_jB_j(t)+v_i
\right)\,\mathrm dt\\
&=
\sum_{i=1}^N v_i
\sum_{A\subset J_i}
\int_{v_A}^{v_A+v_i}
t^{|A|}(1-t)^{N-1-|A|}\,\mathrm dt\\
&=
\sum_{i=1}^N v_i
\int_0^{v_i}
\sum_{A\subset J_i}
(t+v_A)^{|A|}(1-t-v_A)^{N-1-|A|}
\,\mathrm dt.
\end{align*}
The sum in the integrand suggests a shifted version of the binomial theorem. Such an identity is provided by Hurwitz’s multinomial identity  \cite[Eq.~(IV)]{Hurwitz1902}:
\[
\sum_{A\subset J_i}
(t+v_A)^{|A|}(1-t-v_A)^{N-1-|A|}
=
\sum_{A\subset J_i}|A|!\prod_{j\in A}v_j.
\]
It follows that
\begin{align*}
\mathfrak B(v)
&=
\sum_{i=1}^N v_i^2
\sum_{A\subset J_i}|A|!\prod_{j\in A}v_j\\
&=
\sum_{i=1}^N v_i^2
\int_0^\infty e^{-s}
\prod_{j\ne i}(1+v_js)\,\mathrm ds\\
&=
\int_0^\infty\left(1-\prod_{i=1}^N e^{-v_i s}(1+v_i s)\right)\,\frac{\mathrm ds}{s^2}.\qedhere
\end{align*}
\end{proof}

\begin{theorem}[$\ell^2$-estimate for the Bernoulli functional]
\label{thm:bernoulli-l2}
For every $v\in\mathbb R^N$, it holds that
\[
\mathfrak B(v)\le \sqrt{\frac\pi2}\,\|v\|_{\ell^2}.
\]
Moreover, the constant $\sqrt{\pi/2}$ is optimal over $N$.
\end{theorem}

\begin{proof}
The case $v=0$ is immediate.
By Theorem~\ref{thm:sign-alignment}, we may assume that $v_i\ge 0$. Then 
\[
\mathfrak B(v)
\le
\int_0^\infty\frac1{s^2}\left(1-e^{-s^2\|v\|_{\ell^2}^2/2}\right)\mathrm ds
=
\sqrt{\frac\pi2}\,\|v\|_{\ell^2}.
\]

The constant $\sqrt{\pi/2}$ is optimal. Indeed, if $v^{(N)}:=(N^{-1/2},\ldots,N^{-1/2})\in\mathbb R^N$, then $\|v^{(N)}\|_{\ell^2}=1$, and dominated convergence yields
\[
\mathfrak B(v^{(N)})
=
\int_0^\infty\frac1{s^2}\left(1-e^{-\sqrt Ns}\left(1+\frac{s}{\sqrt N}\right)^N\right)\mathrm ds
\to
\sqrt{\frac\pi2}.\qedhere
\]
\end{proof}

\subsection{Subquadratic Bernoulli inequalities}
\label{subsec:subquadratic-bernoulli}

The same integral representation of $\mathfrak B$ allows us to
determine, for each $1\le q\le2$, the optimal dimension-free
constant $C_q$ in the inequality
\[
\mathfrak B(v)\le C_q\|v\|_{\ell^q}.
\]
We have shown that $C_2=\sqrt{\pi/2}$.

\begin{theorem}[Structure of maximizers]
\label{thm:structure-of-maximizers}
For every $1\le q< 2$,
if $v\in\mathbb R^N$ maximizes the functional $\mathfrak B$ on the unit $\ell^q$-sphere, then there exists $1\le k\le N$ such that, up to a global sign and a permutation of coordinates,
\[
v=\frac{1}{k^{1/q}}(\underbrace{1,\dots,1}_{k},0,\dots,0).
\]
\end{theorem}

To prove Theorem~\ref{thm:structure-of-maximizers}, we show that all nonzero coordinates of a maximizer must be equal. Otherwise, the lemma below provides a perturbation that preserves the $\ell^q$-norm while increasing $\mathfrak{B}$, contradicting maximality.

\begin{lemma}
\label{lem:two-coordinate-variation}
Let $x>y>0$ and let $v\in\mathbb R^N$ with $v_i\ge 0$.
For every $1\le q<2$, the functional
\[
\Phi(\varepsilon):=\mathfrak B((x^q+\varepsilon)^{1/q},(y^q-\varepsilon)^{1/q},v)
\]
cannot attain a local maximum at $\varepsilon=0$.
\end{lemma}

\begin{proof}
For sufficiently small $|\varepsilon|$, we write $x_\varepsilon:=(x^q+\varepsilon)^{1/q}$, $y_\varepsilon:=(y^q-\varepsilon)^{1/q}$, and
\[
W_\varepsilon(s):=
e^{-(x_\varepsilon+y_\varepsilon)s}
\prod_{i=1}^N e^{-v_i s}(1+v_i s).
\]
Exponential decay, uniform for small $|\varepsilon|$, permits differentiation under the integral and gives
\begin{align*}
q\Phi'(\varepsilon)
&=
(x_\varepsilon^{2-q}-y_\varepsilon^{2-q})\int_0^\infty W_\varepsilon(s)\,\mathrm ds\\
&\qquad
-(x_\varepsilon y_\varepsilon^{2-q}-x_\varepsilon^{2-q}y_\varepsilon)\int_0^\infty sW_\varepsilon(s)\,\mathrm ds.
\end{align*}
The endpoint case $q=1$ is immediate: $\Phi'(0)>0$, so $\varepsilon=0$ is not a local maximum of $\Phi$.

Assume henceforth that $1<q<2$. Since $x_\varepsilon y_\varepsilon^{2-q}-x_\varepsilon^{2-q}y_\varepsilon>0$, the sign of $\Phi'(\varepsilon)$ is determined by
\[
\sgn \Phi'(\varepsilon)
=
\sgn\left(\frac{(x_\varepsilon + y_\varepsilon)(x_\varepsilon^{2-q}-y_\varepsilon^{2-q})}{x_\varepsilon y_\varepsilon^{2-q}-x_\varepsilon^{2-q} y_\varepsilon}
-
(x_\varepsilon + y_\varepsilon)\frac{\int_0^\infty sW_\varepsilon(s)\,\mathrm ds}{\int_0^\infty W_\varepsilon(s)\,\mathrm ds}
\right).
\]
On the one hand, to show that the first term is strictly increasing in $\varepsilon$, set $r:=x_\varepsilon/y_\varepsilon>1$.
Writing $u=\tfrac12\log r>0$, we have
\[
\frac{(r+1)(r^{2-q}-1)}{r-r^{2-q}}
=\frac{\sinh((3-q)u)}{\sinh((q-1)u)}-1.
\]
The logarithmic derivative of the ratio is positive because $3-q>q-1>0$ and $z\coth z$ is strictly increasing: its derivative is $(\sinh z\cosh z-z)/\sinh^2z>0$.
Since $\mathrm dr/\mathrm d\varepsilon>0$, the chain rule yields
\[
\frac{\mathrm d}{\mathrm d\varepsilon}\frac{(x_\varepsilon+y_\varepsilon)(x_\varepsilon^{2-q}-y_\varepsilon^{2-q})}
{x_\varepsilon y_\varepsilon^{2-q}-x_\varepsilon^{2-q}y_\varepsilon}
=
\frac{\mathrm d}{\mathrm dr}\frac{(r+1)(r^{2-q}-1)}{r-r^{2-q}}\,\frac{\mathrm d}{\mathrm d\varepsilon}\frac{x_\varepsilon}{y_\varepsilon}
>
0.
\]
On the other hand, to show that the second term is nonincreasing in \(\varepsilon\), we use the logarithmic derivative of \(W_\varepsilon\) and integration by parts to obtain
\[
\frac{W_\varepsilon'(s)}{W_\varepsilon(s)}=-(x_\varepsilon+y_\varepsilon)-\sum_{i=1}^N\frac{v_i^2s}{1+v_i s},
\qquad
\int_0^\infty sW_\varepsilon'(s)\,\mathrm ds=-\int_0^\infty W_\varepsilon(s)\,\mathrm ds.
\]
Normalize $W_\varepsilon(s)\,\mathrm ds$ to a probability measure $\nu_\varepsilon$.
The function $h(s)=\sum_i v_i^2s^2/(1+v_i s)$ is nondecreasing.
Thus $\frac{\mathrm d}{\mathrm d\varepsilon}\mathbb E_{\nu_\varepsilon}h=-(x_\varepsilon'+y_\varepsilon')\operatorname{Cov}_{\nu_\varepsilon}(h,s)\ge0$, since $x_\varepsilon'+y_\varepsilon'=(x_\varepsilon^{1-q}-y_\varepsilon^{1-q})/q<0$.
Consequently,
\[
\frac{\mathrm d}{\mathrm d\varepsilon}\left((x_\varepsilon+y_\varepsilon)\frac{\int_0^\infty sW_\varepsilon(s)\,\mathrm ds}{\int_0^\infty W_\varepsilon(s)\,\mathrm ds}\right)
=
-\frac{\mathrm d}{\mathrm d\varepsilon}\left(\frac{\int_0^\infty\sum_{i=1}^N\frac{v_i^2s^2}{1+v_i s}W_\varepsilon(s)\,\mathrm ds}{\int_0^\infty W_\varepsilon(s)\,\mathrm ds}\right)
\le
0.
\]
Therefore, the expression inside the final $\sgn$ is strictly increasing in $\varepsilon$.
If $0$ were a local maximum of $\Phi$, then $\Phi'(0)=0$, so this expression would vanish at $0$.
Consequently $\Phi'(\varepsilon)<0$ for small $\varepsilon<0$ and $\Phi'(\varepsilon)>0$ for small $\varepsilon>0$, making $0$ a strict local minimum, a contradiction.
\end{proof}

\begin{proof}[Proof of Theorem~\ref{thm:structure-of-maximizers}]
By Theorem~\ref{thm:sign-alignment}, all nonzero coordinates
of $v$ have the same sign.
Replacing $v$ by $-v$ if necessary, we may assume
that they are positive.
It remains to show that they are equal.

Suppose they are not, and, up to a permutation, write $v=(x,y,\widetilde v)$ with $x>y>0$.
Since $(x^q+\varepsilon)+(y^q-\varepsilon)=x^q+y^q$,
the variation
\[
\Phi(\varepsilon)
=
\mathfrak B((x^q+\varepsilon)^{1/q},(y^q-\varepsilon)^{1/q},\widetilde v)
\]
preserves the unit $\ell^q$-sphere.
The maximality of $v$ thus implies that $\varepsilon=0$ is a local maximum of $\Phi$, contradicting Lemma~\ref{lem:two-coordinate-variation}.
\end{proof}

Theorem~\ref{thm:structure-of-maximizers} reduces the determination of $C_q$ to a discrete optimization problem.
Indeed, writing $\mathbf 1_N=(1,\ldots,1)\in\mathbb R^N$ and $S_N:=\mathfrak B(\mathbf 1_N)$, we have
\[
C_q=\sup_{N\ge1}\frac{S_N}{N^{1/q}}.
\]
It therefore remains to determine the optimal $N$ as $q$ varies.

\begin{theorem}[Exact phase transitions]
\label{thm:exact-phase-transitions}
Define the sequence $\{q_N\}_{N=1}^\infty$ by 
\[
\frac{1}{q_N}:=\frac{\log S_{N+1}-\log S_N}{\log(N+1)-\log N}.
\]
Then the following statements hold:
\begin{enumerate}
\item The sequence $\{q_N\}_{N=1}^\infty$ increases strictly to $2$.
\item For every \(q_{N-1}\le q\le q_N\) (with the convention \(q_0:=1\)), it holds that
\[
C_q = \frac{S_N}{N^{1/q}}.
\]
\item Every maximizer has exactly $N$ nonzero coordinates when $q_{N-1}<q<q_N$, and $N$ or $N+1$ nonzero coordinates when $q=q_N$.
\end{enumerate}
\end{theorem}

\begin{proof}
\textbf{Step 1.} To compare the critical values $q_N$, we interpolate $S_N$ on the logarithmic scale by setting
\[
S(x):=\int_0^\infty e^{-s}
\left(1+\frac{s}{e^x}\right)^{e^x-1}\,\mathrm ds,
\]
so that \(S_N=S(\log N)\). The change of variables $1+s/e^x=e^{e^y}$ yields
\[
S(x)=\int_{\mathbb R}e^{H(x,y)}\,\mathrm dy,
\qquad
H(x,y)=x+y-\sum_{k=2}^\infty\frac{e^{x+ky}}{k!}.
\]
The Hessian of $H$ has quadratic form $-\sum_{k\ge2}e^{x+ky}(a+kb)^2/k!$, which is strictly negative for $(a,b)\ne(0,0)$ by the terms $k=2,3$.
The strictness in the Pr\'ekopa--Leindler argument can be seen directly. Fix $x_0\ne x_1$ and $0<\theta<1$, and normalize $\rho_j(y)=e^{H(x_j,y)}/S(x_j)$.
The increasing transport $T$ between these positive smooth densities satisfies $\rho_1(T(y))T'(y)=\rho_0(y)$.
Changing variables by $y\mapsto(1-\theta)y+\theta T(y)$, strict concavity of $H$ gives
\begin{align*}
S((1-\theta)x_0+\theta x_1)
&>S(x_0)^{1-\theta}S(x_1)^\theta
\int_{\mathbb R}\rho_0(y)\frac{(1-\theta)+\theta T'(y)}{T'(y)^\theta}\,\mathrm dy\\
&\ge S(x_0)^{1-\theta}S(x_1)^\theta,
\end{align*}
where the last inequality is the weighted arithmetic--geometric mean inequality.
The first inequality is strict because $(x_0,y)\ne(x_1,T(y))$ for every $y$. Thus $S$ is strictly log-concave.
The secant slopes of $\log S$ at the points $\log N$ therefore decrease strictly, namely
\[
\frac1{q_N}=\frac{\log S_{N+1}-\log S_N}{\log(N+1)-\log N}>\frac1{q_{N+1}}.
\]

\textbf{Step 2.}
By the integral representation for $\mathfrak B$,
\[
\frac{S_N}{\sqrt N}
=\int_0^\infty\frac1{s^2}
\left(
1-e^{-\sqrt Ns}\left(1+\frac{s}{\sqrt N}\right)^N
\right)\mathrm ds.
\]
The integrands increase strictly to \((1-e^{-s^2/2})/s^2\): this follows from the strict decrease of $(u-\log(1+u))/u^2=\int_0^1r/(1+ur)\,\mathrm dr$ for $u>0$.
Monotone convergence consequently gives
\[
\frac{S_N}{\sqrt N}
\uparrow
\int_0^\infty\frac{1-e^{-s^2/2}}{s^2}\,\mathrm ds
=\sqrt{\frac\pi2}.
\]
The strict monotonicity of \(S_N/\sqrt N\) yields \(1/q_N>1/2\), and the strict log-concavity of \(S\) yields, for \(N\ge2\),
\[
\frac1{q_N}
\le
\frac{\log S_N-\log S_1}{\log N}
\to
\frac12.
\]
Combining these bounds shows that \(q_N\to2\), completing the proof of (1).

\textbf{Step 3.} It remains to maximize \(S_N/N^{1/q}\) over \(N\). Observe that
\[
\log\frac{S_{N+1}/(N+1)^{1/q}}{S_N/N^{1/q}}
=
\log\frac{S_{N+1}}{S_N}-\frac1q\log\frac{N+1}{N}
=
\left(\frac1{q_N}-\frac1q\right)\log\frac{N+1}{N}.
\]
The sign of the preceding expression determines the optimal support sizes. Assertions (2) and (3) then follow from Theorem~\ref{thm:structure-of-maximizers}.
\end{proof}

The first few phases are summarized in Table~\ref{tab:first-five-phases}.
Numerically, $q_1\approx1.709511$, $q_2\approx1.758887$, $q_3\approx1.787342$, $q_4\approx1.806858$, and $q_5\approx1.821441$.

\begin{table}[htbp]
\centering
\begingroup
\setlength{\tabcolsep}{10pt}
\renewcommand{\arraystretch}{1.5}
\begin{tabular}{@{}ccc@{}}
\toprule
Range of $q$ & Optimal support size & $C_q$ \\
\midrule
$1\le q\le q_1$ & $1$ & $S_1=1$ \\[0.5em]
$q_1\le q\le q_2$ & $2$ &
$\frac{S_2}{2^{1/q}}=\frac{3}{2\cdot2^{1/q}}$ \\[0.5em]
$q_2\le q\le q_3$ & $3$ &
$\frac{S_3}{3^{1/q}}=\frac{17}{9\cdot3^{1/q}}$ \\[0.5em]
$q_3\le q\le q_4$ & $4$ &
$\frac{S_4}{4^{1/q}}=\frac{71}{32\cdot4^{1/q}}$ \\[0.5em]
$q_4\le q\le q_5$ & $5$ &
$\frac{S_5}{5^{1/q}}=\frac{1569}{625\cdot5^{1/q}}$ \\
\bottomrule
\end{tabular}
\endgroup
\vspace{0.75em}
\caption{The first five phases of the dimension-free optimization problem. At each transition point $q=q_k$, both support sizes $k$ and $k+1$ are optimal.}
\label{tab:first-five-phases}
\end{table}

\subsection{Sharp Poisson inequalities}
For $\mathfrak P_T$, Cauchy--Schwarz gives
\[
\mathfrak P_T(v)
\le
\|v\|_{\ell^2}\int_0^T\frac{\mathrm ds}{\sqrt s}=2\sqrt T\|v\|_{\ell^2}.
\]
To obtain the sharp constant, we realize $P_i(s)$ in the definition of $\mathfrak P_T$ as the values at time $s$ of independent Poisson processes of unit rate.

\begin{proposition}
For every $v\in\mathbb R^N$ and $T>0$, it holds that
\[
\mathfrak P_T(v)=2\mathbb E\left[-\inf_{0\le s\le T}\sum_{i=1}^Nv_i(P_i(s)-s)\right].
\]
\end{proposition}

\begin{proof}
Since the increments of $X_s:=\sum_{i=1}^Nv_i(P_i(s)-s)$ on the $n$-grid are independent and identically distributed,
Spitzer's combinatorial identity \cite[Eq.~(3.1)]{Spitzer1956}, applied to \(-X\) on the grid \(\{iT/n:0\le i\le n\}\), gives
\[
\mathbb E\left[-\min_{0\le i\le n}X_{iT/n}\right]
=
\sum_{i=1}^{n}\frac{\mathbb E[(-X_{iT/n})^+]}{i}
=
\frac T{2n}\sum_{i=1}^{n}
\frac{\mathbb E|X_{iT/n}|}{iT/n}.
\]
Since \(s\mapsto \mathbb E|X_s|/s\) is continuous on \((0,T]\) and bounded by \(\|v\|_{\ell^2}s^{-1/2}\), the right-endpoint Riemann sums converge to the improper integral, giving
\[
\lim_{n\to\infty}
\frac T{2n}\sum_{i=1}^{n}
\frac{\mathbb E|X_{iT/n}|}{iT/n}
=
\frac12\int_0^T\mathbb E|X_s|\,\frac{\mathrm ds}s
=
\frac12\mathfrak P_T(v).
\]
Meanwhile, right-continuity of $X$ and monotone convergence over successively refined dyadic partitions give
\[
\lim_{r\to\infty}
\mathbb E\left[-\min_{0\le i\le2^r}X_{iT/2^r}\right]
=
\mathbb E\left[-\inf_{0\le s\le T}X_s\right].\qedhere
\]
\end{proof}

\begin{theorem}
For every $v\in\mathbb R^N$ and $T>0$, it holds that
\[
\mathfrak P_T(v)\le \sqrt{\frac{8T}{\pi}}\|v\|_{\ell^2}.
\]
Moreover, the constant $\sqrt{8T/\pi}$ is optimal over $N$.
\end{theorem}

\begin{proof}
The case $v=0$ is immediate. By Theorem~\ref{thm:sign-alignment} and homogeneity, we may assume that $v_i\ge0$ and $\sum_{i=1}^Nv_i^2=1$.
Introduce the Gaussian smoothing
\[
\varphi(s,x):=\mathbb E|x+\sqrt{T-s}Z|,
\qquad
Z\sim N(0,1),
\]
so that $\varphi$ is smooth for $x\in\mathbb R$ and $s<T$, and satisfies $\varphi_x(s,0)=0$.
Next, define
\[
R_s:=X_s-\inf_{0\le u\le s}X_u,
\qquad
X_s:=\sum_{i=1}^Nv_i(P_i(s)-s),
\]
so that $\mathfrak P_T(v)=2\mathbb ER_T$.
For fixed $0<\sigma<T$, It\^o's formula for the finite-variation semimartingale $\varphi(s,R_s)$ gives\begin{align*}
\varphi(\sigma,R_\sigma)-\varphi(0,0)
&=
\int_0^\sigma\left[\varphi_s(s,R_{s-})
-\sum_{i=1}^Nv_i\varphi_x(s,R_{s-})\right]\,\mathrm ds\\
&\qquad
+\sum_{i=1}^N\int_0^\sigma
\left[\varphi(s,R_{s-}+v_i)-\varphi(s,R_{s-})\right]\,\mathrm dP_i(s).
\end{align*}
Taking expectations and applying the Poisson compensation formula to the jump integrals yields
\begin{align*}
\mathbb E\varphi(\sigma,R_\sigma)-\sqrt{\frac{2T}{\pi}}
&=\mathbb E\int_0^\sigma\left[\varphi_s(s,R_{s-})-\sum_{i=1}^N v_i\varphi_x(s,R_{s-})\right]\,\mathrm ds\\
&\qquad
+\sum_{i=1}^N \mathbb E\int_0^\sigma\left[\varphi(s,R_{s-}+v_i)-\varphi(s,R_{s-})\right]\,\mathrm dP_i(s)\\
&=\int_0^\sigma\mathbb E\left[\varphi_s(s,R_{s-})-\sum_{i=1}^N v_i\varphi_x(s,R_{s-})\right.\\
&\qquad \left.+\sum_{i=1}^N \varphi(s,R_{s-}+v_i)-\sum_{i=1}^N\varphi(s,R_{s-})\right]\,\mathrm ds\\
&\le 0.
\end{align*}
To justify the last inequality, note that $\varphi_{xx}(s,x)=-2\varphi_s(s,x)$, and that $\varphi_{xx}(s,\cdot)$ is nonincreasing on $[0,\infty)$. Thus, Taylor's formula in the $x$-variable gives
\[
\varphi(s,x+v_i)-\varphi(s,x)
\le
v_i\varphi_x(s,x)+\frac{v_i^2}{2}\varphi_{xx}(s,x),
\qquad x\ge0.
\]
Almost surely, $R$ is continuous at time $T$, and $\varphi(\sigma,R_\sigma)\to R_T$ as $\sigma\uparrow T$.
Fatou's lemma therefore yields
\[\mathfrak P_T(v)
\le
2\liminf_{\sigma\uparrow T}\mathbb E\varphi(\sigma,R_\sigma)
\le
\sqrt{\frac{8T}\pi}.
\]

The constant $\sqrt{8T/\pi}$ is optimal.
Indeed, if $v^{(N)}:=(N^{-1/2},\ldots,N^{-1/2})\in\mathbb R^N$, then $\|v^{(N)}\|_{\ell^2}=1$, and the central limit theorem, together with uniform integrability and dominated convergence, yields
\[
\mathfrak P_T(v^{(N)})
=
\int_0^T\frac1s\,\mathbb E\left|\frac1{\sqrt N}\sum_{i=1}^N(P_i(s)-s)
\right|\,\mathrm ds
\to
\sqrt{\frac{8T}\pi}.\qedhere
\]
\end{proof}

\section{Beckmann Poincar\'e inequalities}

Combining the Beckmann duality theorem with the explicit Malliavin weights in Section 3 and the sharp integral inequalities in Section 5 yields Beckmann Poincaré inequalities on Hamming cubes and product Poisson spaces, as well as a Beckmann pseudo-Poincaré inequality on finite abelian Cayley graphs. 

Throughout this section, \(C_q\) denotes the optimal Bernoulli constant determined in Theorems 5.2 and 5.5. We first specify the normalization in each upper bound, then prove the stated sharpness properties of the Beckmann inequalities.

\begin{corollary}[Hamming cubes]
\label{cor:hamming-cubes}
For \(1\le q\le2\), the Beckmann Poincaré inequality
\[
\|g-\mu(g)\|_{L^1(\mu)}
\le
2^{2/q-1}C_q\,\mathsf B_q(g),\qquad g:\{-1,1\}^N\to \mathbb R
\]
holds on all Hamming cubes $\{-1,1\}^N$.
Moreover, the constant \(2^{2/q-1}C_q\) is optimal uniformly over \(N\).
\end{corollary}

\begin{proof}
Here $b(x,y)=\mu(x)/2$, so $|v|_q(x)=2^{-2/q}\|v(x,\cdot)\|_{\ell^q}$.
For $T=\infty$, Example~\ref{ex:cube-weights} therefore bounds the optimal constant by $\frac12\sup_{|v|_q(x)\le1}\mathfrak B(v)\le2^{2/q-1}C_q$.

Fix $N,d\ge 1$. Let $v_{N,d}:\{-1,1\}^{N+d}\to \mathbb R^{N+d}$ be the vector field defined by
\[
v_{N,d}(x)
:=
\begin{cases}
2(\mathbf 1_N,\mathbf 0_d)&x=\mathbf 1_{N+d},\\
\mathbf 0_{N+d}&x\not=\mathbf 1_{N+d},
\end{cases}
\]
and let $g_{N,d}$ be the mean-zero solution of $\Delta g_{N,d}=\operatorname{div}v_{N,d}$.
Then, by the definition of \(\mathsf B_q\), we have
\[
2^{N+d}\mathsf B_q(g_{N,d})
\le 
2^{1-2/q}N^{1/q}.
\]

To obtain a lower bound for \(\|g_{N,d}\|_{L^1(\mu)}\), 
let $X_t^{\mathbf 1_{N+d}}$ and $H_t(\mathbf 1_{N+d})$ be as in Example~\ref{ex:cube-weights}, and
define
\[
f_{N,d}(x)=-\sgn\left(\sum_{i=1}^N x_i-\frac{N}{d}\sum_{i=N+1}^{N+d}x_i\right).
\]
Using duality and the Malliavin weights, we obtain
\[
\begin{aligned}
2^{N+d}\langle g_{N,d},f_{N,d}\rangle_\mu
&=
2^{N+d}\int_0^\infty\langle v_{N,d},\nabla P_tf_{N,d}\rangle_b\,\mathrm dt\\
&=
\frac12
\int_0^\infty
\sum_{i=1}^N
(\nabla P_tf_{N,d}(\mathbf1_{N+d}))_i\,\mathrm dt\\
&=
\frac12
\int_0^\infty
\mathbb E\left[
f_{N,d}(X_t^{\mathbf1_{N+d}})
\sum_{i=1}^N (H_t(\mathbf1_{N+d}))_{i}
\right]\,\mathrm dt.
\end{aligned}
\]
Since \(\|f_{N,d}\|_\infty\le1\), the law of large numbers and dominated convergence yield, for each fixed \(N\),
\[
\liminf_{d\to\infty}
2^{N+d}\|g_{N,d}\|_{L^1(\mu)}
\ge
\lim_{d\to\infty}
2^{N+d}\langle g_{N,d},f_{N,d}\rangle_\mu
=
\mathfrak B(\mathbf1_N).
\]
Therefore, for \(1\le q\le 2\), the constant \(2^{2/q-1}C_q\) is optimal
\[
\sup_{N\ge1}\liminf_{d\to\infty}
\frac{\|g_{N,d}\|_{L^1(\mu)}}{\mathsf B_q(g_{N,d})}
\ge
2^{2/q-1}C_q.\qedhere
\]
\end{proof}

\begin{corollary}[Product Poisson spaces]
\label{cor:product-poisson-spaces}
For $1\le q\le 2$, the Beckmann Poincar\'e inequality
\[
\|g-\mu(g)\|_{L^1(\mu)}
\le
2^{1/q}C_q\mathsf B_q(g),
\qquad
g\in D(\Delta)
\]
holds on all product Poisson spaces \(\mathbb N_0^N\).
Moreover, for \(q=2\), the constant \(\sqrt\pi\) is optimal, even in dimension two.
\end{corollary}

\begin{proof}
Writing $v_i^\pm=v(x,x\pm e_i)$, with $v_i^-=0$ if $x_i=0$, we have
\[
|v|_q(x)^q=\frac12\left(\sum_i|v_i^+|^q+\sum_i x_i|v_i^-|^q\right).
\]
Conditional on $M_i=P_i(1)$, the bridge $P_i(s)-sP_i(1)$ is a sum of $M_i$ centered Bernoulli variables, while $B_i^x(s)-x_is$ is a sum of $x_i$ such variables, all independent.
At $T=\infty$, the two integrals in Example~\ref{ex:poisson-weights} coincide and cancel the factor $1/2$; the conditional Bernoulli bound and Jensen's inequality give
\[
C_q\mathbb E\left(\sum_i M_i|v_i^+|^q+\sum_i x_i|v_i^-|^q\right)^{1/q}
\le2^{1/q}C_q|v|_q(x),
\]
since $\mathbb EM_i=1$. This proves the upper bound.

Fix $d\ge 1$ and put $x_d:=(d,d)$.
Let $v_d:\mathbb N_0^2\times\mathbb N_0^2\to\mathbb R$ be the vector field defined by
\[
v_d(x,y)
:=
\begin{cases}
d^{-1/2}&(x,y)=(x_d,x_d-e_1),\\
-d^{-1/2}&(x,y)=(x_d,x_d-e_2),\\
0&\text{otherwise},
\end{cases}
\]
and let $g_d$ be the mean-zero solution of $\Delta g_d=\operatorname{div}v_d$.
Then, by the definition of $\mathsf B_2$, we have
\[
\frac{\mathsf B_2(g_d)}{\mu(x_d)}\le1.
\]

To obtain a lower bound for \(\|g_d\|_{L^1(\mu)}\), let \(X_t^{x_d}\) and \(H_t\) be as in Example~\ref{ex:poisson-weights}, and define
\[
f(x_1,x_2):=\operatorname{sgn}(x_2-x_1).
\]
Using duality and the weight formula, we obtain
\[
\begin{aligned}
\frac{\langle g_d,f\rangle_\mu}{\mu(x_d)}
&=
\frac1{\mu(x_d)}
\int_0^\infty
\langle v_d,\nabla P_tf\rangle_b\,\mathrm dt\\
&=
\frac{\sqrt d}{2}\int_0^\infty
\sum_{i=1}^2(-1)^{i+1}
\nabla P_tf(x_d,x_d-e_i)\,\mathrm dt\\
&=
\frac{\sqrt d}{2}\int_0^\infty
\mathbb E\left[
f(X_t^{x_d})
\sum_{i=1}^2(-1)^{i+1}H_t(x_d,x_d-e_i)
\right]\,\mathrm dt.
\end{aligned}
\]
Put $D_d(s)=B_1^{x_d}(s)-B_2^{x_d}(s)$ and $R(s)=P_2(s)-P_1(s)$.
Since $\mathrm dt=\mathrm ds/(1-s)$, the last integral equals
\[
\int_0^1\frac{\mathbb E[D_d(s)\sgn(D_d(s)+R(s))]}{2\sqrt d\,s(1-s)}\,\mathrm ds.
\]
For fixed $0<s<1$, $D_d(s)/\sqrt d$ converges to $N(0,2s(1-s))$, whereas $R(s)/\sqrt d\to0$ in $L^2$; the second-moment bound gives uniform integrability.
The integrand thus converges to $1/\sqrt{\pi s(1-s)}$ and is bounded in absolute value by $1/\sqrt{2s(1-s)}$, since $\mathbb E|D_d(s)|\le\sqrt{2ds(1-s)}$.
Dominated convergence gives $\pi^{-1/2}\int_0^1[s(1-s)]^{-1/2}\,\mathrm ds=\sqrt\pi$; hence, since $\|f\|_\infty\le1$,
\[
\liminf_{d\to\infty}
\frac{\|g_d\|_{L^1(\mu)}}{\mu(x_d)}
\ge
\lim_{d\to\infty}
\frac{\langle g_d,f\rangle_\mu}{\mu(x_d)}=
\sqrt\pi.
\]
Therefore, for $q=2$, the constant $\sqrt\pi$ is optimal
\[
\liminf_{d\to\infty}
\frac{\|g_d\|_{L^1(\mu)}}{\mathsf B_2(g_d)}
\ge\sqrt\pi.\qedhere
\]
\end{proof}

\begin{corollary}[Abelian Cayley graphs]
\label{cor:abelian-cayley-graphs}
For every \(T>0\), the Beckmann pseudo-Poincaré inequality
\[
\|g-P_Tg\|_{L^1(\mu)}\le\sqrt{\frac{4T}\pi}\,\mathsf B_2(g),
\qquad
g:\Gamma\to \mathbb R,
\]
holds on all finite abelian Cayley graphs $G=\operatorname{Cay}(\Gamma,S)$.
Moreover, the constant \(\sqrt{4T/\pi}\) is optimal, even among cyclic Cayley graphs.
\end{corollary}

\begin{proof}
Here $|v|_2(x)=\frac12\|v(x,\cdot)\|_{\ell^2}$.
Example~\ref{ex:abelian-cayley-weights} and the sharp Poisson estimate bound the optimal constant by $\frac12\sup_{|v|_2(x)\le1}\mathfrak P_{T/2}(v)\le\sqrt{4T/\pi}$.

Fix $d\ge 1$. Let $G_d=\operatorname{Cay}(\mathbb Z_{2d^2+1},\{\pm 1,\dots,\pm d\})$, and let $v_d$ be the vector field defined by
\[
v_d(x,y):=\frac{\sqrt2y}{\left(\sum_{i=1}^d i^2\right)^{1/2}},
\qquad
x=0,
\qquad 
1\le |y|\le d,
\]
with $v_d(x,y)=0$ on all other edges, and let $g_d$ be the mean-zero solution of $\Delta g_d=\operatorname{div}v_d$.
Here \(y\) is identified with its integer representative in \(\{-d^2,\ldots,d^2\}\).
Put $A_d=\sum_{i=1}^d i^2$, so that $\mathsf B_2(g_d)\le\mu(0)$, and set $f_d(z)=\sgn(\widetilde z)$ for the integer representative $\widetilde z\in[-d^2,d^2]$.
For independent $P_i^\pm(t/2)\sim\operatorname{Pois}(t/2)$, let $Y_{d,t}=\sum_{i=1}^d i(P_i^+(t/2)-P_i^-(t/2))$, whose variance is $tA_d$.
Duality and Example~\ref{ex:abelian-cayley-weights} give
\[
\frac{\|g_d-P_Tg_d\|_{L^1(\mu)}}{\mu(0)}
\ge\frac{\langle g_d-P_Tg_d,f_d\rangle_\mu}{\mu(0)}
=\int_0^T\frac{\mathbb E[Y_{d,t}\sgn(\widetilde{Y_{d,t}})]}{\sqrt{2A_d}\,t}\,\mathrm dt.
\]
For fixed $t>0$, the characteristic-function expansion
\[
\log\mathbb E e^{iuY_{d,t}/\sqrt{tA_d}}
=t\sum_{i=1}^d\left(\cos\frac{ui}{\sqrt{tA_d}}-1\right)
=-\frac{u^2}{2}+O_{t,u}(d^{-1})
\]
and uniform integrability give $\mathbb E|Y_{d,t}|/\sqrt{tA_d}\to\sqrt{2/\pi}$.
The cyclic wrap error satisfies
\[
\left|\mathbb E[Y_{d,t}\sgn(\widetilde{Y_{d,t}})]-\mathbb E|Y_{d,t}|\right|
\le2\mathbb E[|Y_{d,t}|\mathbf1_{\{|Y_{d,t}|>d^2\}}]
\le\frac{2tA_d}{d^2}=o(\sqrt{A_d}).
\]
Finally, Cauchy--Schwarz bounds the absolute integrand by $(2t)^{-1/2}$.
Dominated convergence gives $\int_0^T(\pi t)^{-1/2}\,\mathrm dt=\sqrt{4T/\pi}$, and hence
\[
\liminf_{d\to\infty}
\frac{\|g_d-P_Tg_d\|_{L^1(\mu)}}{\mathsf B_2(g_d)}
\ge
\sqrt{\frac{4T}{\pi}}.\qedhere
\]
\end{proof}

We conclude by verifying the sharpness of the classical cube gradient inequality at $q=2$ stated in the introduction.
Set $S_N(x)=N^{-1/2}\sum_i x_i$ and $g_{N,a}=h_a(S_N)$, where $h_a(u)=\tanh(au)$ and $a>0$.
Symmetry gives $\mu(g_{N,a})=0$, while Taylor's formula gives
$\bigl||\nabla_2g_{N,a}|-h_a'(S_N)\bigr|\le\|h_a''\|_\infty/\sqrt N$.
For $Z\sim N(0,1)$, the central limit theorem therefore yields
\[
\lim_{a\to\infty}\lim_{N\to\infty}
\frac{\|g_{N,a}\|_{L^1(\mu)}}{\mathbb E|\nabla_2g_{N,a}|}
=\lim_{a\to\infty}\frac{\mathbb E|h_a(Z)|}{\mathbb E h_a'(Z)}
=\sqrt{\frac\pi2},
\]
since the numerator tends to $1$ and the denominator to $2/\sqrt{2\pi}$.

\section*{Acknowledgements}
This work grew out of the course ``Foundations of High-Dimensional Probability'', taught by the first author at Shandong University in spring 2026, and a seminar at University of Science and Technology of China. We thank Ruizhe Zhang for helpful discussions and the  participants in the seminar at USTC for their questions and comments.

\section*{Use of AI tools}
The authors used ChatGPT and Codex to explore proof strategies, develop and check arguments, and revise the manuscript. In July 2026, we obtained a combinatorial identity for centered Bernoulli sums with nonnegative coefficients. Ruizhe Zhang provided ideas from analytic number theory and Fourier methods that motivated our use of Fourier positivity for coefficients of mixed signs. Building on these ideas, we carried out multiple rounds of exploration and calculation with the assistance of ChatGPT, producing a first draft in early August 2026. We subsequently continued to refine and extend the paper, also using Kimi during this process. The research objectives, the course of the investigation, and the selection and formulation of the results were determined entirely by the authors, based on their own interests and judgements. The authors can fully explain the rationale behind every step throughout the research.

We subsequently became aware of an online preprint addressing the same sharp Hamming-cube inequality. The overlap with our first draft is substantial: it includes the Fourier sign-alignment argument through sinc-squared distributions, reflection inequalities under convolution and a positive sine kernel; the exact integral representation for nonnegative coefficients; the Gaussian comparison giving the sharp quadratic constant; and the finite-dimensional values corresponding to equal coefficients. These are central components of the quadratic Bernoulli argument.

\end{document}